\documentclass[referee]{gji}
\usepackage{timet,color}
\usepackage[urlcolor=blue,citecolor=black,linkcolor=black]{hyperref}
\usepackage{apalike}
\let\bibhang\relax
\usepackage{natbib}
\usepackage{graphicx}
\usepackage{color}
\usepackage{listings}
\usepackage{amsmath}

\usepackage{amsthm}
\usepackage{amssymb}
\usepackage{amsbsy}
\usepackage{float}
\usepackage{wasysym}
\usepackage{setspace}
\usepackage{bm}
\usepackage{algorithm}
\usepackage{algpseudocode}
\usepackage{caption}
\usepackage{subcaption}
\usepackage{comment}
      \newcommand{\Ja}{J_{\alpha}}
      \newcommand{\tJa}{\tilde{J}_{\alpha}}	
      \newcommand{\bx}{{\bf x}}
      
      \newcommand{\aw}{w_{\alpha}}

      \newtheorem{result}{Result}
      \newtheorem{definition}{Definition}

\begin{document}

\title[Solution of an Inverse Problem by Extended Inversion]
  {Solution of an Acoustic Transmission Inverse Problem by Extended Inversion}
\author[William W. Symes et al.]
  {William W. Symes$^1$,  Huiyi Chen$^2$ and Susan E. Minkoff$^2$ \\
  $^1$ PO Box 43, Orcas, WA 98280 USA \\
  $^2$ Department of Mathematical Sciences, FO 35, 
  University of Texas at Dallas, Richardson, TX 75080 USA
  }

\pubyear{2021}

\maketitle
\begin{summary}
Study of a simple single-trace transmission example shows how an
extended source formulation of full-waveform inversion can produce an
optimization problem without spurious local minima (“cycle skipping”),
hence efficiently solvable via Newton-like local optimization
methods. The data consist of a single trace extracted from a causal
pressure field, propagating in a homogeneous fluid according to linear
acoustics, and recorded at a given distance from a transient point
energy source. The source intensity (“wavelet”) is presumed
quasi-impulsive, with zero energy for time lags greater than a
specified maximum lag. The inverse problem is: from the recorded
trace, recover both the sound velocity or slowness and source wavelet
with specified support, so that the data is fit with prescribed RMS
relative error. The least-squares objective function has multiple
large residual minimizers. The extended inverse problem permits source
energy to spread in time, and replaces the maximum lag constraint by a
weighted quadratic penalty. A companion paper shows that for proper
choice of weight operator, any stationary point of the extended
objective produces a good approximation of the global minimizer of the
least squares objective, with slowness error bounded by a multiple of
the maximum lag and the assumed noise level. This paper summarizes
the theory developed in the companion paper and presents numerical
experiments demonstrating the accuracy of the predictions in concrete
instances. We also show how to dynamically adjust the penalty scale
during iterative optimization to improve the accuracy of the slowness
estimate.

\end{summary}


\section{Introduction}
Model-based parameter estimation by least squares data fitting is a
widely used and successful technique for data analysis in science and
engineering \cite[]{Bertero:98,Vogel:02}, and particularly in the
geosciences \cite[]{Parker:94,Tarantola:05}. In its application to
seismology, least-squares data fitting has come to be known as Full
Waveform Inversion (FWI) and is now well-established as a useful tool
for probing the Earth's subsurface
\cite[]{VirieuxOperto:09,Fichtner:10,Schuster:17}. However, FWI encounters a serious
practical obstacle, so-called ``cycle-skipping''. This term signifies
the tendency of the mean-square error function to exhibit approximate stationary points far from a global optimum.
%(we
%show an explicit example of this phenomenon later in this paper). 
Because of the
typical dimensions of field FWI problems, only local gradient-based
minimization algorithms are computationally practical. These methods find
approximate stationary points (or local minima), so may stagnate at suboptimal and geologically
uninformative Earth models. Local descent methods avoid suboptimal
stagnation only if initial models are already quite close to optimal,
in the sense of predicting the arrival times of seismic events to
within a small multiple of the dominant data wavelength
\cite[]{GauTarVir:86,VirieuxOperto:09,Plessix:10}.

This paper examines one of a number of proposed alternatives to least
squares FWI, based on an extension of the model-to-data mapping
(``extended inversion''), in the context of 
a very simple acoustic transmission inverse problem: the sound
velocity in a homogeneous fluid is to be determined from a single recording of a
pressure wave emitted at a known
location around a known time, and propagated over a known distance. In a companion paper \cite[]{Symes:22a}, we show precisely
how descent-based FWI fails to solve this problem absent a very good
initial model estimate and precisely why extended inversion succeeds, using the same type
of optimization, and from an essentially arbitrary initial model
estimate. The companion paper gives detailed estimates for the accuracy of the sound
velocity determination in terms of the time extent of the energy
excitation and the signal-to-noise ratio of the data. In this paper we describe these theoretical results succinctly and illustrate the mathematical conclusions with relevant numerical
examples.

Model extension is an old idea in
seismic data processing 
\cite[]{geoprosp:2008}. In recent years, many extended inversion algorithms have been proposed as cures for cycle-skipping
\cite[]{SymesCar:91,Plessix:00a,Symes:09,LuoSava:11,BiondiAlmomin:SEG12,LeeuwenHerrmannWRI:13,LiuSymesLi:14,Warner:14,ChaurisGP:14,Warner:16,LeeuwenHerrmann:16,Herve2017,HouSymes:Geo18,Aghmiryetal:20,MetivierBrossier:SEG20}. The approach studied in this paper uses {\em
  source extension}: the energy source component of the experimental model is permitted to
have more (or less constrained) parameters than the experimental
design suggests. \cite{HuangNammourSymesDollizal:SEG19} give an
overview and taxonomy of source extension methods. 

This work 
%also explores an extended inversion application, but 
differs in a crucial respect from the works cited above. In all of those papers the authors use numerical tests with synthetic and/or field data to justify their conclusions about how the methods avoid cycle-skipping. Of course, numerical examples provide retrospective evidence, which in itself cannot guarantee that another test will not produce a substantially different result. In fact, \cite{Symes:2020} showed that one of these approaches, believed for several years to be immune to cycle-skipping, is in some cases no better in that respect than FWI.  

The theoretical results presented in the companion paper \cite[]{Symes:22a} provide the missing guarantee. These results, reviewed in the next section, explain circumstances under which the extended inversion algorithm developed below will succeed in producing an accurate model estimate starting with an arbitrary initial guess. Our numerical examples illustrate essential features of this theory. Some of the error bounds are sharp, and we give examples that underline such conclusions. In other cases the theory gives sufficient conditions  
which may not be necessary under practically interesting limitations on the data beyond the basic assumptions of the theory. We illustrate such opportunities for further theoretical development as well.

Extended inversion is not the only alternative to straightforward
least-squares data fitting that may overcome cycle-skipping. For
example, evidence has recently emerged that the Wasserstein metric
arising in the theory of optimal transport may provide a measure
of error between model-predicted and observed seismic data that is less
oscillatory than the mean-square error
\cite[]{EngquistYang:GEO18,Metivier:GEO18}. Very recently, \cite{MahankaliYang:21} have developed theoretical results analogous to ours, showing convexity of the Wasserstein-2 error function in several simple examples.

\begin{comment}
The single-trace transmission problem studied in this paper is a
drastically simplified cartoon of actual seismic prospecting. It seems
worth careful study nonetheless for two reasons: it exhibits the
phenomenon of cycle-skipping, and its simplicity allows the properties
of both FWI and extended inversion approaches to be established, both theoretically and by numerical example.
\end{comment}

In the next (``Theory'') section, we describe the conceptual
framework of this inverse problem, review the related literature, and
state the essence of our main theoretical results as established in \cite{Symes:22a}. 
The theory is illustrated in the ``Numerical Experiments'' Section that
follows.

\section{Theory}

In this section we
describe the theoretical results which we will illustrate with numerical experiments in the next section.
\subsection{Background and Problem Description}

The physical setting of our study is constant-density acoustic wave propagation in a {3D} homogeneous medium in which the
(single) mechanical parameter of interest used to describe the Earth is the speed of sound, or for our purposes the reciprocal wave velocity or
slowness, $m>0$.

The source of acoustic
energy is localized at a point in space and 
radiates uniformly in all
directions with time-dependent intensity (``wavelet'') $w(t)$, thereby generating a pressure
field, $p(\bx,t)$. Both the 
pressure field and the wavelet vanish at large
negative times. The data consists of a recording of the pressure field at a single receiver located a distance, 
$r$, from the source, over the time interval $[t_{\rm min}$,
$t_{\rm max}]$. The pressure field
solves the wave equation \cite[]{Frie:58}:
\begin{eqnarray}
  \label{eqn:awe}
  \left(m^2\frac{\partial^2 }{\partial t^2} - \nabla^2\right) p(\bx,t) &=&
                                                                         w(t)\delta(\bx-\bx_s) \nonumber\\
  p(\bx,t)&=&0, \,\,\,\,t\ll 0.
\end{eqnarray}
The solution $p$ of the problem \ref{eqn:awe} is given by the well-known expression
(\cite{CourHil:62}, Chapter VI, section 12, equation 47):
\begin{equation}
  \label{eqn:homsol}
  p(\bx,t) = \frac{1}{4\pi |\bx-\bx_s|}w\left(t-m|\bx-\bx_s|\right).
\end{equation}
The predicted pressure at the receiver point, $\bx_r$, is
\begin{equation}
\label{eqn:mod}
p(\bx_r,t) = \frac{1}{4\pi r}w\left(t-mr\right)  = F[m]w(t).
\end{equation}
The mapping $F[m]$ so defined (ignoring the amplitude factor
$1/(4\pi r)$) is an $m$-dependent time shift, which is the basis of many illustrations of
the cycle-skipping phenomenon (for example, see \cite{VirieuxOperto:09},
Figure 7). It is unsurprising, therefore, that an analysis of cycle-skipping
can be rooted in the simplified setting we describe. 
A naive version of the inverse problem investigated in this paper
would be: given recorded data
$d,$ 
find $m$ and $w$ 
so that
$F[m]w \approx d$. 

Clearly this problem statement is not interesting as it
stands, since 
fitting
the data does not constrain the slowness $m$ at all: for any $m$, the choice $w(t) = 4\pi rd(t+mr)$
yields $F[m]w = d$ exactly.
To be useful, the problem statement must be augmented with a
constraint on the 
model $(m,w)$.  One natural constraint is to assume
that $w$ is non-zero only for small times: a maximum time lag (``support radius'') $\lambda>0$ exists so that $w(t)=0$ for $|t|>\lambda$. 

An assumption that $w$ is non-zero only for small $|t|$ may be justified as
follows. Seismic sources may act
over considerable time intervals. An example from exploration
seismology is the airgun,
which produces an oscillating pressure pulse that dies away
slowly. It is commonplace to estimate this pulse, usually by recording
the emitted pressure field at a position where it should nominally be isolated from other
signal, then deconvolve it
from the data by safeguarded Fourier division or other means, with
appropriate compensation for the relative positions of the pulse recording and the zone of interest within the earth.
This so-called signature deconvolution
(\cite{Ziolkowski:GEO82,SheriffGeldart:1995},
section 9.5 and \cite{Yil:01}, Chapter 2)
results in modified data
corresponding to a source equal to signature deconvolution of the
pulse estimate itself. The deconvolved source pulse approximates an
impulse ($\delta(t)$). It cannot actually be an impulse, due to the finite
frequency nature of seismic data,
but its amplitudes are typically small outside of a 
much smaller time interval than is the case for the original pulse
estimate. 

This reasoning leads 
us to incorporate specification of a
maximum time lag or support radius
$\lambda$ into the statement of the inverse problem. To state the problem precisely, we use the standard notation $\|u\|^2$ for the square integral (``$L^2$ norm squared") of any function $u$ of time:
\begin{equation}
    \label{eqn:normdef}
    \|u\|^2 = \int_{-\infty}^{\infty}dt\,u(t)^2.
\end{equation}

\noindent {\bf Inverse Problem:}
{\em given data 
    $d$, assumed data noise level $\epsilon \in
  [0,1)$, extremal slownesses $0< m_{min} \leq m_{max}$ and support radius $\lambda 
  > 0$, find $m \in [m_{min}, m_{max}]$ and $w$ vanishing for $|t| > \lambda$, so that
\begin{equation}
  \label{eqn:probstat0}  \|F[m]w-d\| \le \epsilon\|d\|.
\end{equation}
}

\subsection{Full Waveform Inversion}
Despite its simplicity, this problem exhibits the 
pathology characteristic of field-scale FWI when 
formulated as a least-squares optimization problem.

Given a data trace $d$
and maximum lag $\lambda>0$, we introduce three versions of the least squares (FWI) objective function
for the single-trace transmission problem:

\begin{definition}
\label{thm:fwibas}
The {\em basic FWI objective} function $e$ of slowness $m$ and wavelet $w$ is
\begin{equation}
  \label{eqn:redms}
   e[m,w;d]=\frac{1}{2}\frac{\|F[m]w-d\|^2}{\|d\|^2}
\end{equation}
\end{definition}

\begin{definition}
\label{thm:fwired}
The {\em reduced FWI objective} $\tilde{e}$ of slowness $m$ is 
\begin{equation}
\label{eqn:vpmfwi}
\tilde{e}[m;d] = \min_w e[m,w;d],
\end{equation}
 the minimum being taken over all $w$ with $w(t)=0$ for $|t|>\lambda$.
\end{definition}

\noindent
\begin{definition}
\label{thm:fwires}
Given a wavelet $w_*$ with $w_*(t)=0$ for $|t|>\lambda$, the {\em restricted FWI objective} with $w=w_*$ is the function of $m$ given by $e[m,w_*;d]$.
\end{definition}

\begin{result}
\label{thm:result1} 
Suppose that $m_*>0$, $\lambda>0$, $w_*(t)=0$ for $|t|>\lambda$, and $d=F[m_*]w_*$. If $m \ne m_*$, then $\tilde{e}[m;d]$ and $e[m,w_*;d]$ are both $>0$. If $|m-m_*| \ge 2\lambda/r$, then
\[
\tilde{e}[m;d]= e[m,0;d] = 1/2
\]
and
\[
e[m,w_*;d]=1.
\]
\end{result}

That is, both the reduced FWI objective $\tilde{e}$ (minimizer of $e$ over $w$ for each $m$) and the restricted FWI objective $e[\cdot,w_*;d]$ (using the same wavelet $w_*$ as used to generate the noise-free data) have entire intervals of local minimizers for $m$ sufficiently far from the slowness $m_*$ used to generate the data, which is also the global minimizer of both. Under the conditions we have posed, the minimum in equation \eqref{eqn:vpmfwi} is attained. Therefore any local minimizer of $\tilde{e}$ corresponds to at least one local minimizer of $e$. Consequently $e$ also has a continuum of local minimizers. These local optima are far from the global minimizer $m_*$ (or global minmizer $(m_*,w_*)$ in the case of $e$) and have 50\% and 100\% error levels respectively, for noise-free data. Local optimization of any of these objectives therefore fails to solve the Inverse Problem stated above for any $\epsilon<1/2$, unless the initial estimate of $m$ is within $2\lambda/r$ of the global minimizer $m_*$. 
This phenomenon is what is meant by ``cycle-skipping''.

The observation, that (locally) minimizing $\tilde{e}$ over $m$ corresponds to minimizing $e$ over $(m,w)$, suggests a {\em nested} approach: inner minimization over $w$ (to compute $\tilde{e}$) and outer minimization (of $\tilde{e})$ over $m$. This nested algorithm is is an example of the {\em Variable Projection Method} (or ``VPM''), in the nomenclature introduced by \cite{GolubPereyra:73,GolubPereyra:03}. Several authors have observed that VPM may be more computationally efficient for FWI problems (including source estimation) than direct application of local optimization algorithms to $e$. In fact, the sensitivity of $e$ to perturbations in source parameters (like $w$) tends to be very different from its sensitivity to perturbations in kinematic medium parameters (such as $m$). The nested iteration segregates these sensitivities and reduces ill-conditioning and resulting computational inefficiency \cite[]{Aravkin:12,LiRickettAbubakar:13,YinHuang:16}.
VPM is the main computational device used in our examples. It also plays a central role in the theory developed in the companion paper.

\subsection{Extended Source Inversion}

``Extension'', as used in this paper, means provision of additional degrees of freedom in the domain of the modeling operator (such as $F$) with the intent of generating more paths avoiding local minima and leading to a solution of the Inverse Problem. 
The extension discussed here consists in dropping the
support constraint on $w$, thus considerably enlarging the space of feasible models. In fact, as noted above, the extended model space is then so large that perfect data fit is possible with any $m$.
To compensate for the resulting indeterminacy in the
estimation of $m$, we add a 
quadratic 
penalty term to the 
least squares 
objective,
$e,$
defined by
a penalty operator $A$
and penalty weight $\alpha$.

\begin{definition}
\label{thm:esibas}
The Extended Source Inversion (``ESI'') objective function $\Ja$ is defined by
\begin{eqnarray}
  g[w] & = & \frac{1}{2}\|Aw\|^2/\|d\|^2, \label{eqn:pen} \\
  \Ja[m, w;d] &=& e[m,w;d] + \alpha^2 g[w] \nonumber\\
  &=& \frac{1}{2}(\|F[m]w-d\|^2 + \alpha^2 \|Aw\|^2)/\|d\|^2.
\label{eqn:penerr}
\end{eqnarray}
\end{definition}
\noindent
The normalizations by the data norm $\|d\|$ in the equations
\eqref{eqn:redms}, \eqref{eqn:pen}, and \eqref{eqn:penerr} make the values
of $e$, $\alpha^2 g$, and $\Ja$ nondimensional and simplify
the interpretation of
numerical examples.

The penalty operator $A$ must 
penalize 
nonzero values of $w$ for large $|t|$ since we search for a wavelet that is only active in a short time interval centered at $t=0$. 
One simple choice for $A$ is 
\begin{equation}
\label{eqn:annihilator}  
(Aw)(t)=tw(t).
\end{equation}
We will assume $A$ is of this form throughout the paper. 
The penalty weight $\alpha$ also must be chosen, and we will discuss the choice of $\alpha$ in the next subsection.

While minimization of the ESI objective $\Ja$ might be tackled directly - by
alternating minimizations over $m$ and $w$ or by computing updates for $m$ and $w$ simultaneously - we have already noted that $e$, a component of $\Ja$, has dramatically different
sensitivity to $m$ versus $w$, and therefore so does $\Ja$. Therefore we employ the Variable Projection Method described above, with an inner minimization of $\Ja$ over $w$, followed by an outer minimization over $m$. 
Any minimizer $w$ of $\Ja[m,\cdot; d]$ must solve the {\em normal equation}
\begin{equation}
  \label{eqn:norm}
  (F[m]^TF[m]+\alpha^2A^TA)w= F[m]^Td. 
\end{equation}
With the definitions of $F$ and $A$ given above, the normal equation has a unique solution. Therefore $\Ja[m,\cdot;d]$ has a unique minimizer $w = w_{\alpha}[m;d]$. 

\begin{definition}
\label{thm:esired}
The {\em reduced ESI objective} $\tJa$ is
\begin{equation}
  \label{eqn:redexp}
  \tJa[m;d] = \inf_w \Ja[m,w;d] = \Ja[m,w_{\alpha}[m;d];d].
\end{equation}
\end{definition}

\begin{result}
\label{thm:result2}
Let $d_*=F[m_*]w_*$ be noise-free data with {\em target slowness} $m_* > 0$ and {\em target wavelet} $w_*(t) = 0$ for $|t| > \lambda$. Let $d=d_*+n$ be noisy data, with noise trace $n$. 
Define the {\em noise-to-signal ratio} $\eta$ by $\eta = \|n\|/\|d_*\|$. If 
$\alpha>0$ and $\eta < \frac{\sqrt{5}-1}{2}$, 
then any stationary point $m$ of $\tJa[\cdot;d]$ satisfies
  \begin{equation}
    \label{eqn:merr}
     |m-m_*| \le (1+f(\eta))\frac{\lambda}{r},
  \end{equation}
  where $f(\eta) = \frac{2\eta(1+\eta)}{1-\eta(1+\eta)} = 2\eta + O(\eta^2)$.
\end{result}

In particular, if the data is noise-free, that is $\eta=0$, then the error between any stationary point of the reduced ESI objective $\tJa$ and the target slowness $m_*$ is at most the maximum lag $\lambda$ of the target wavelet $w_*$ divided by the source-receiver offset $r$. Note the dramatic contrast with the least-squares estimation (minimization of $e$ as in \ref{eqn:redms}). As indicated in Result \ref{thm:result1}, the reduced FWI objective has stationary points at arbitrarily large differences from the target slowness.

No general bounds can hold for much
larger noise levels than indicated in Result \ref{thm:result2}. For instance, $n=-F[m_*]w_*$ ($\eta=1$) yields $d=0$, for which any $m$ is a stationary point.
However, numerical exploration suggests that stronger bounds might hold given other constraints on data error. A natural example is uniformly distributed random noise
filtered to have the same spectrum as the source. Numerical Experiment 5 below illustrates this case, for which the only stationary point of the reduced ESI objective is a quite precise estimator of the target slowness. 

Unless the data is noise-free, the
estimated wavelet $\aw[m;d]$ (solution of the normal equation \eqref{eqn:norm}) corresponding to a stationary point $m$ of $\tJa[\cdot;d]$ need not vanish for $|t|>\lambda$. To construct a solution
of the Inverse Problem, we must modify
the wavelet $w_\alpha[m,d]$. 

\begin{result}
    \label{thm:result3} Define \[
{\bf 1}_{[-\lambda, \lambda]}(t) = 
 \begin{cases}
      1, & |t| \leq \lambda \\
      0, & \text{else}.
 \end{cases} 
 \]
Suppose that $\mu>0$, $d = F[m_*]w_* + n$ with $w_*(t) = 0$ for $|t| > \mu$, and $m$ is a stationary point of $\tJa[\cdot;d]$. Then the pair $(m, {\bf 1}_{[-\lambda, \lambda]} w_\alpha[m,d])$ solves the Inverse Problem \eqref{eqn:probstat0}, that is
\[
   e[m,{\bf 1}_{[-\lambda,\lambda]}w;d] \le \frac{1}{2} \epsilon^2,
\]
provided that
\begin{align}
    \label{eqn:lambbd}
    \lambda \ge (2+f(\eta))\mu,\\ 
    \label{eqn:epsbd}
    \epsilon \ge \frac{(8 \pi r \alpha \mu)^2}{1+(8 \pi r \alpha \mu)^2} + \eta.
\end{align}
\end{result}
Thus the minimization of the extended objective $J_\alpha$ provides a solution of the Inverse Problem with constraints on the achievable performance and bounds on the errors.

\subsection{The discrepancy algorithm}

So far in this story, the selection of the penalty weight $\alpha$ has played a relatively minor role. On the other hand, in practical calculations selection of $\alpha$ strongly influences the convergence of iterative methods and the quality of the results.
Concentration of the wavelet near $t=0$ clearly favors larger $\alpha$. In
fact the reduced penalty term $g[m,\aw[m;d];d]$ is a decreasing function of $\alpha$ for any $m$. As $g$ measures the dispersion of the wavelet away from $t=0$, larger $\alpha$ is to be preferred, all else being equal.
However all else is not equal: the data error $e[m,\aw[m;d];d]$ is an increasing function of $\alpha$, so $\alpha$ is constrained by data fit. The choice of $\alpha$ affects the reduced ESI objective $\tJa[m;d]$ and, therefore, the estimation of the nonlinear variable $m$ as well:
\begin{result}
    \label{thm:result4}
Assume once again that $d_*=F[m_*]w_*$, $d=d_*+n$ with $w_*(t)=0$ for $|t|>\lambda$ and noise-to-signal ratio $\eta=\|n\|/\|d_*\|$, and suppose that $m$ is a stationary point of $\tJa[\cdot;d]$. Then
\begin{equation}
    \label{eqn:mdiff}
    |m -m_*| \le \frac{\lambda}{r} + \frac{\eta}{\alpha}(...),
\end{equation}
in which the elided quantity is polynomial in $\alpha$ with positive constant term.
\end{result}
\noindent
That is, stationary points of $\tJa$ become more accurate approximations of the target slowness $m_*$ as $\alpha$ increases, at least for small $\alpha$.

These observations motivate an algorithm for selection of $\alpha$,
based on a version of the {\em discrepancy principle}
\cite[]{EnglHankeNeubauer,Hanke:17,FuSymes2017discrepancy}, adapted to the setting of this paper, namely $m$ and $\alpha$ should solve the 
\vspace{.1in}

\noindent {\bf Discrepancy Problem:}
{\em
  given data $d \in D$ and a range of minimum and maximum allowable errors $0 < e_- < e_+$,  find the slowness $m$ and the scalar $\alpha$ so that
  \begin{itemize}
  \item[(i) ]$m$ is a stationary point of the reduced objective function $\tJa[\cdot;d]$, and
  \item[(ii) ]$e_- < e[m,\aw[m;d];d] < e_+$.
  \end{itemize}
  }
\noindent  
(``Discrepancy'' is used here as a synonym for ``data misfit'', as measured by $e$.)

We use an alternating, or coordinate search, {\it discrepancy algorithm} for
solution of this problem, combining a local optimization
algorithm for updating $m$, and an algorithm for updating
$\alpha$. A first version of the discrepancy algorithm appeared in \cite{FuSymes2017discrepancy}. 

Note that the mean square error, $e,$ lies in the interval $(e_-,e_+)$ at a solution of this problem. Use of an interval, rather than a single target error level, accomplishes two objectives:
\begin{itemize}
    \item it is consistent with the general lack of precise knowledge of data error in most applications;
    \item it permits a local optimization algorithm to make several small updates of $m$ before an update of $\alpha$ is required.
\end{itemize}
In general, lack of an obvious rule for selection of an initial $\alpha$ obstructs the use of the discrepancy principle for penalty formulations of many inverse problems. A common solution is a more or less arbitrary initial selection followed by repeated objective evaluation at a geometric sequence, or a bisection series, of weights until one is found for which the required quantity lies in the specified interval. 

The extended inverse problem studied here has a remarkable property that eliminates this obstruction: a useful initial $\alpha$ value for the problem studied here is simply $\alpha=0$. Since the modeling operator $F[m]$ is surjective for any $m>0$, as noted earlier, this choice yields $e[m,w_0[m;d];d]=0$ for any initial $m$ (which may therefore be chosen arbitrarily). Thus the upper bound on $e$ in condition (ii) of the Discrepancy Problem  is trivially satisfied.

To get the data misfit into the prescribed range (satisfy the lower bound in condition (ii)) and to update it after subsequent updates of $m$, we use the 

\vspace{.1in}
\noindent {\bf Basic Alpha Update Rule:} 
{\em
given an initial value $\alpha_0$ of $\alpha$, generate a sequence $\{\alpha_k: k=0,1,2,...\}$  by the recursion
 \begin{equation}
\label{eqn:alphasecant}
\alpha_{k+1}= \left(\alpha_k^2 + \frac{e_{+}-e[m,w_{\alpha_k}[m;d];d]}{2g[w_{\alpha_k}[m;d]]} \right)^{1/2},
\end{equation}
}
first suggested by \cite{FuSymes2017discrepancy}.

\cite{Symes:21a} shows (Appendix A) that $e$ is increasing in $\alpha$ and that the sequence of $\alpha$ produced by the update rule \ref{eqn:alphasecant} generates an increasing sequence of $e$ in the interval $[0,e_+]$ with $e_+$ as its only accumulation point, hence attaining the target interval $[e_-,e_+]$ in a finite (and estimable) number of steps. Convergence can be accelerated by a number of devices. However in the examples reported in the next section, we have used the basic rule \ref{eqn:alphasecant}.

Having updated $\alpha$, the algorithm updates $m$ by one or more steps of a local optimization algorithm. In the examples presented below, we have used Brent's method \cite[]{Brent:TCJ71} to search for a zero of the gradient.

After the $\alpha$ update cycle, the error bounds are satisfied, that is,
$e[m,\aw[m;d];d] \in (e_-,e_+)$, but the $m$ update is likely to reduce $e$, as it is a
summand in the definition of $\tJa$. If $e$ is reduced below
$e_-$ then control passes to the $\alpha$ update. If an approximate local minimizer is detected, then the algorithm terminates if the error bounds are satisfied, or passes to the $\alpha$ update if not. Thus $\alpha$ is steadily increased, resulting in increasing accuracy of the $m$ estimate, while the data misfit is kept within the prescribed range.
The next section recounts several uses of this discrepancy algorithm, illustrating its behaviour in detail.

\subsection{Solving the Inverse Problem}

Our goal is to solve the Inverse Problem stated at inequality \eqref{eqn:probstat0}. That is, to determine a source $w$ with prescribed support $\subset [-\lambda, \lambda]$ and the slowness $m$ in a prescribed interval $[m_{min}, m_{max}]$, so that the RMS data error $||F[m]w-d||$ is less than a prescribed multiple of $||d||$. Merely finding a stationary point of the reduced ESI objective function $\tJa$ does not accomplish this goal, even though it yields a good slowness estimate under many circumstances, both because the wavelet $\aw[m;d]$ so obtained does not in general vanish for $|t|>\lambda$, and because the stationarity condition does not constrain data error.

We propose an algorithm for solving the Inverse Problem via ESI that addresses both of these shortcomings, consisting of two steps. In the first step we use the discrepancy algorithm, described in the previous subsection, to adjust the penalty weight $\alpha$ dynamically, in the course of minimization of the reduced ESI objective $\tJa$ (thus estimating the wavelet $w=\aw[m;d]$ at the same time). The discrepancy algorithm maintains an upper bound on data error $e$ throughout the minimization of $\tJa$, and simultaneously maximizes $\alpha$, thus minimizing wavelet energy dispersion away from $t=0$. The second step of the algorithm, truncation of the wavelet to the prescribed maximum lag $|t|\le \lambda$ as in Result \ref{thm:result3}, therefore modifies the wavelet as little as possible, and in turn increases the data error as little as possible.

The reasoning just stated is partly heuristic; the theory developed in the companion paper \cite[]{Symes:22a} does not directly address the performance of this algorithm. Result \ref{thm:result3} states sufficient conditions on the data error and maximum lag under which a stationary point of $\Ja$ will solve the Inverse Problem but does not really explain the role of $\alpha$ or of the discrepancy algorithm for controlling it. Instead, in this paper we include a pair of examples, at the end of the next section, which suggest that our approach has some promise. 

\section{Numerical Examples}

This section presents numerical examples that illustrate the theoretical results described in the previous section. The examples are divided into two sets. The first set (Experiments 1-5) explores the location of stationary points for both FWI and ESI. The results of these experiments conform to the predictions of the theoretical Results \ref{thm:result1} and \ref{thm:result2}. The second set (Experiments 6 and 7) present an implementation of the discrepancy algorithm and show how it can be used in conjunction with ESI to solve the Inverse Problem using the algorithm sketched at the end of the previous section. These last two examples illustrate Results \ref{thm:result3} and \ref{thm:result4}.

In each numerical experiment, as in the theory, waves are generated by a single source, and data is recorded at a single receiver. The source and receiver are located $1$ km apart ($r = |x_r - x_s| = 1$ km). Unless noted otherwise, the target wavelet $w_*(t)$ is a the truncation of a zero phase Ricker wavelet 
with peak frequency $40$ Hz, truncated to have support radius $\mu = 0.025$ s (see Figure~\ref{fig:figure_9}). The target slowness is $m_* = 0.4$ s/km. We choose a noise trace $n$, and add it to the noise-free data predicted by $m_*$ and $w_*$, to obtain data traces $d = F[m_*]w_*(t) + n(t)= \frac{1}{4\pi r} w_*(t-m_*r) + n(t)$, recorded over the interval $[0.25,0.65]$ s.

Explicit expressions for the reduced objective function $\tJa$ and its gradient $\frac{d}{dm} \tJa$ are:
\begin{equation}
  \label{eqn:expjgen}
  \tJa[m;d] = \frac{1}{2\|d\|^2}
\int_{t_{\rm min}}^{t_{\rm max}}\,dt\,(4\pi r \alpha (t-mr))^2(1+(4\pi r \alpha 
(t-mr))^2)^{-1}d(t)^2,
\end{equation}

\begin{equation}
  \label{eqn:dexpjgen}
 \frac{d}{dm}\tJa[m;d] =  -\frac{(4 \pi r \alpha)^2}{\|d\|^2}
 \int_{t_{\rm min}}^{t_{\rm max}} \,dt \, 
  (t-mr)(1+(4\pi r \alpha 
  (t-mr))^2)^{-2}d(t)^2. 
\end{equation}

These two expressions are based on Definition \ref{thm:esired} and the Green's function in equation \eqref{eqn:mod}. A detailed derivation can be found in \cite[]{Symes:22a}. Here we approximate these integrals by the trapezoidal rule with $[t_{\rm min},t_{\rm max}] = [0.25, 0.65]$ s.  Since the frequency of the Ricker wavelet is $40$ Hz, its wavelength is $0.025$ s. We pick a step size of $dt = 0.001$ resulting in 25 grid points per wavelength.

\subsection{Stationary Points}
In the five experiments discussed in this subsection, we explore the locations of stationary points for restricted and reduced FWI and reduced ESI objective functions and their dependence on data noise-to-signal level $\eta$ and properties of the target wavelet $w_*$. In all cases 
the penalty weight $\alpha=1$. 

\noindent
\textbf{Experiment 1: Noise Free Data}\\ In this experiment the data is noise-free (noise-to-signal ratio $\eta$ is equal to $0$). Figure~\ref{fig:nf_data} displays the data $d = F[m_*]w_*(t)$. Figure~\ref{fig:figure_1} shows the restricted FWI objective function (Definition \ref{thm:fwires}) using the target wavelet $w_*$, in blue and the reduced ESI objective function (Definition \ref{thm:esired}) in red, as functions of slowness. We see that the restricted FWI objective function has infinitely many local minima, in the regions [0.25, 0.35] and [0.45, 0.63] s/km, consistent with Result 1 with $\lambda=0.025$ s, along with several other stationary points nearer the target slowness. %$m=m_*$ s/km. 
This plot is similar to many examples in the literature illustrating the tendency of FWI to cycle-skip, thus rendering local optimization ineffective. In contrast, the reduced ESI objective function exhibits only one local minimum, at the target slowness, consistent with Result \ref{thm:result2}.

\begin{figure}
\centering
\includegraphics[width=0.5\textwidth]{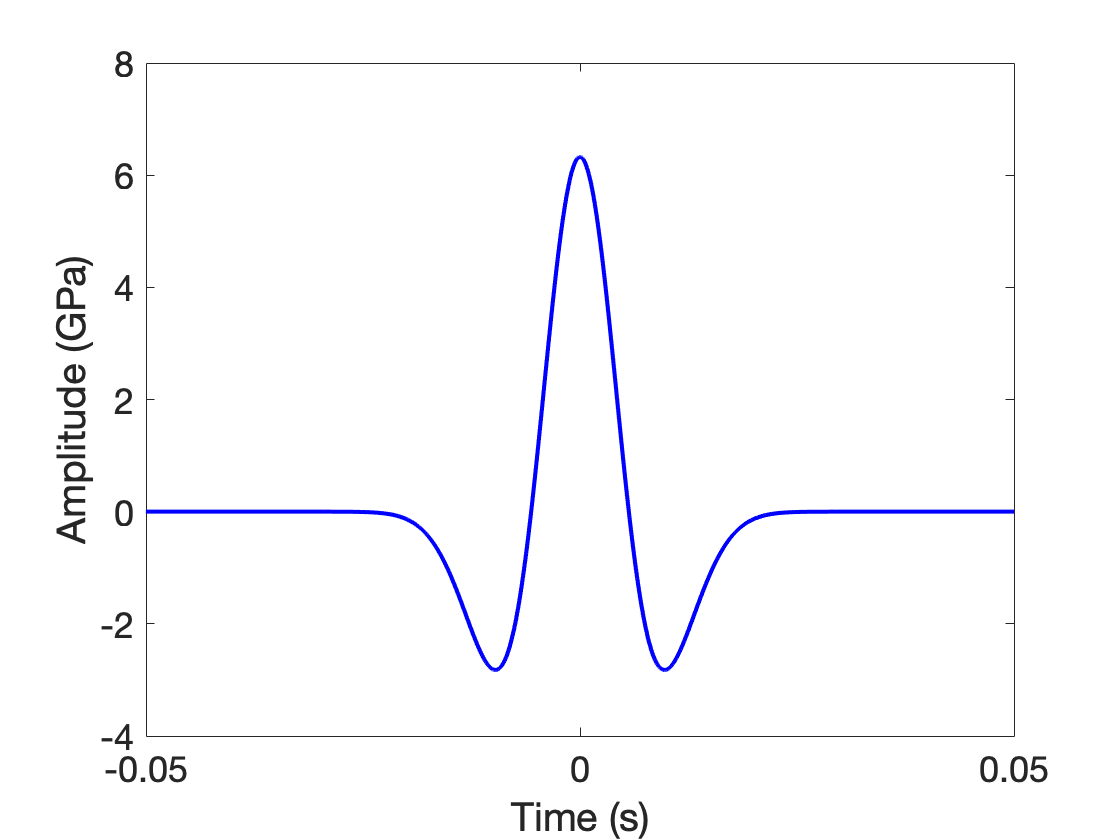}
\caption{$40$ Hz Ricker wavelet centered at $t=0$.}
\label{fig:figure_9}
\end{figure}

\begin{figure}
\centering
\includegraphics[width=0.5\textwidth]{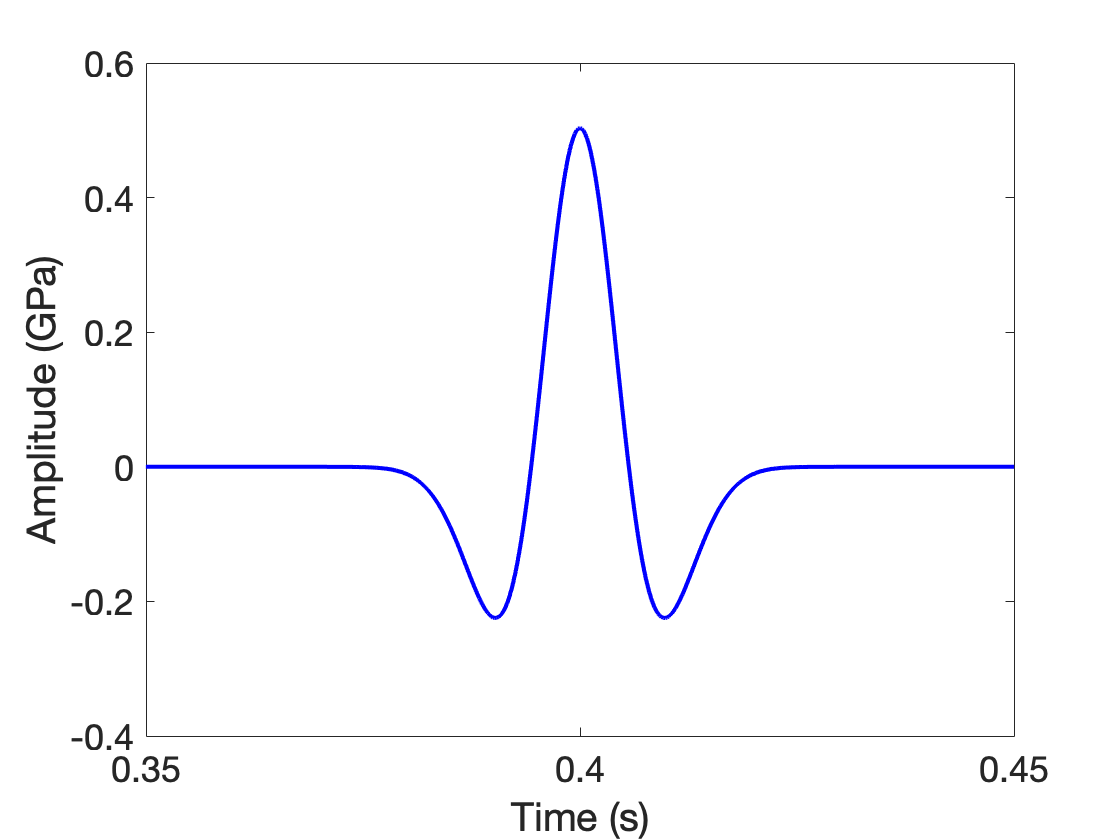}
\caption{Experiment 1: data.}
\label{fig:nf_data}
\end{figure}

\begin{figure}
\centering
\includegraphics[width=0.5\textwidth]{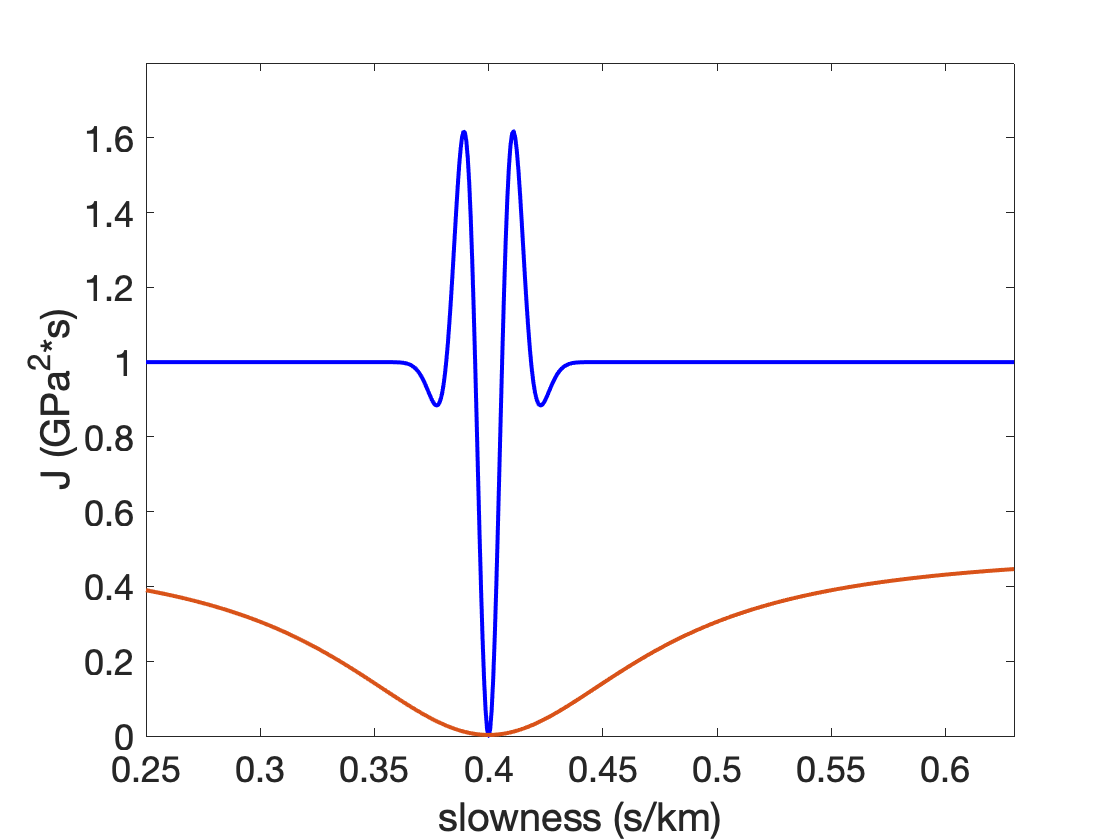}
\caption{Experiment 1: the restricted FWI (blue curve) and reduced ESI (red curve) objective functions versus slowness.}
\label{fig:figure_1}
\end{figure}

%%%%%%%%%%%
%%%%%%%%%%%
\noindent
\textbf{Experiments 2a \& 2b: Noise Free Data, Shifted Wavelet}
These examples show that the bounds mentioned in Result \ref{thm:result2} are sharp and that the reduced (rather than restricted) FWI objective defined in Definition \ref{thm:fwired} behaves as indicated in Result \ref{thm:result1}.

\begin{figure}
\centering
\includegraphics[width=0.5\textwidth]{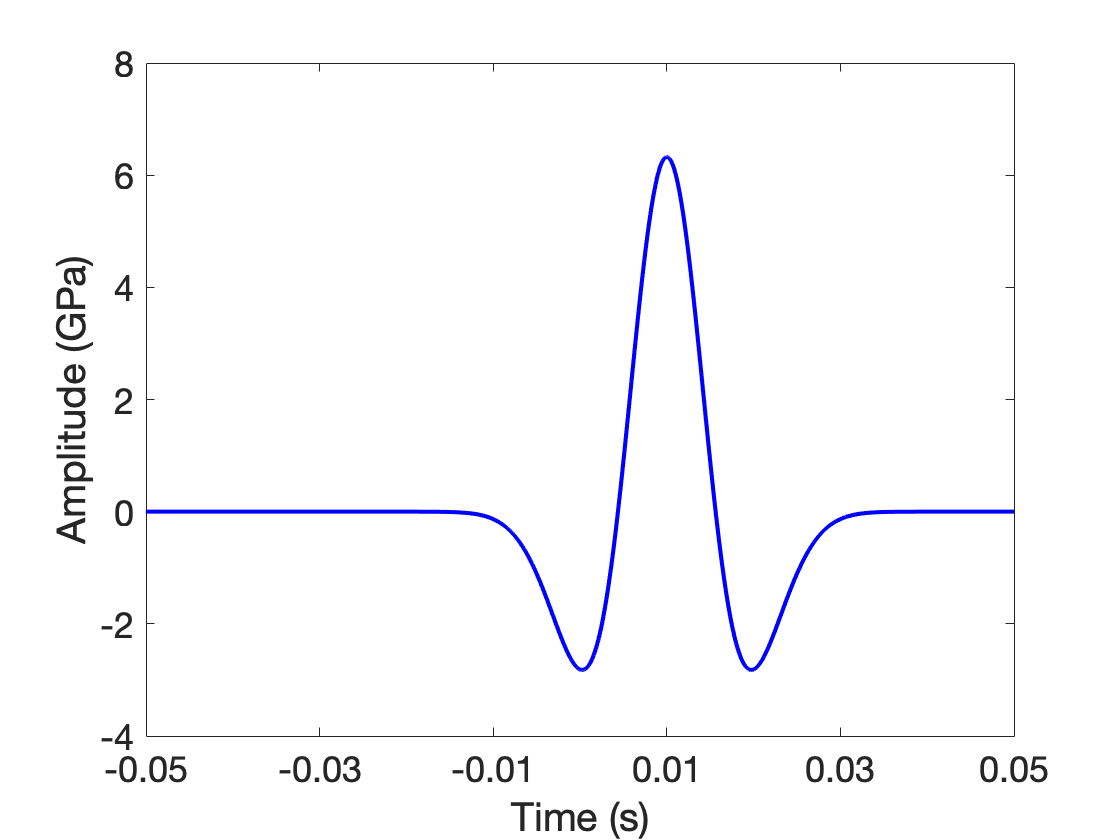}
\caption{Target wavelet for Experiment 2a with peak amplitude at $t=0.01$ s.}
\label{fig:diff_target_wavelet}
\end{figure}

\begin{figure}
\centering
\includegraphics[width=0.5\textwidth]{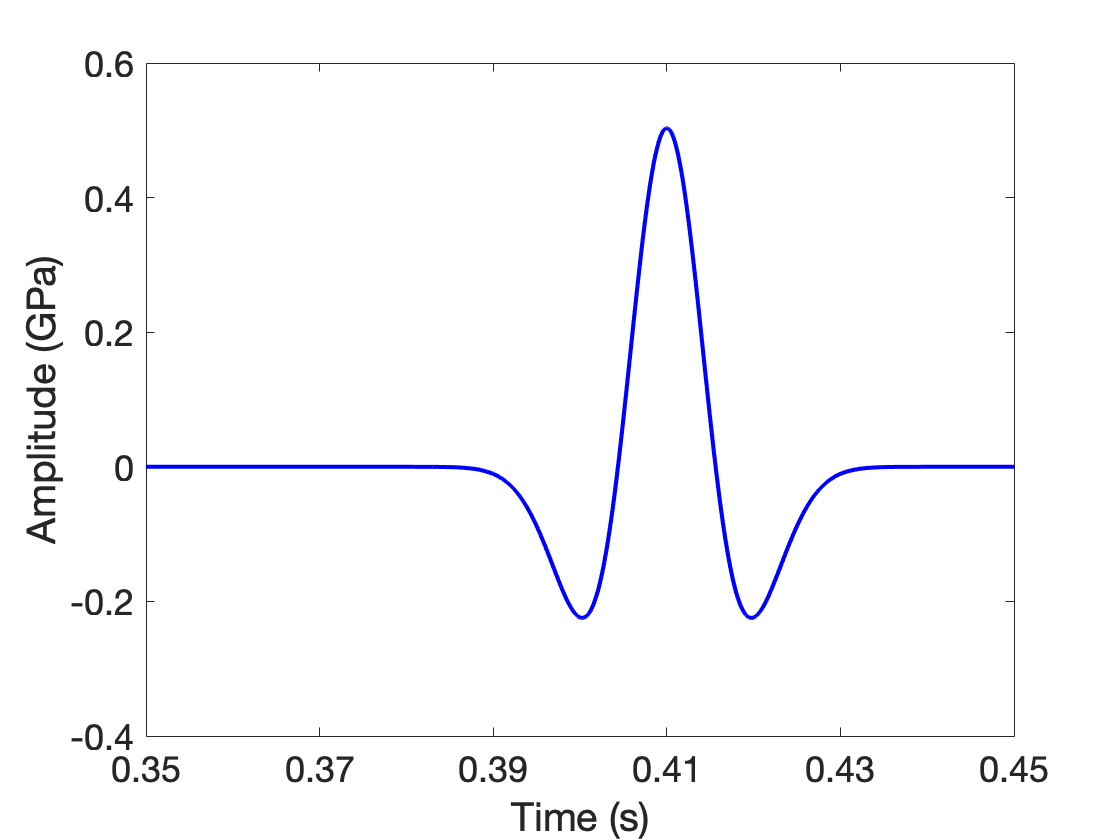}
\caption{Noise-free data used in Experiment 2a.}
\label{fig:data_1}
\end{figure}

\begin{figure}
\centering
\includegraphics[width=0.5\textwidth]{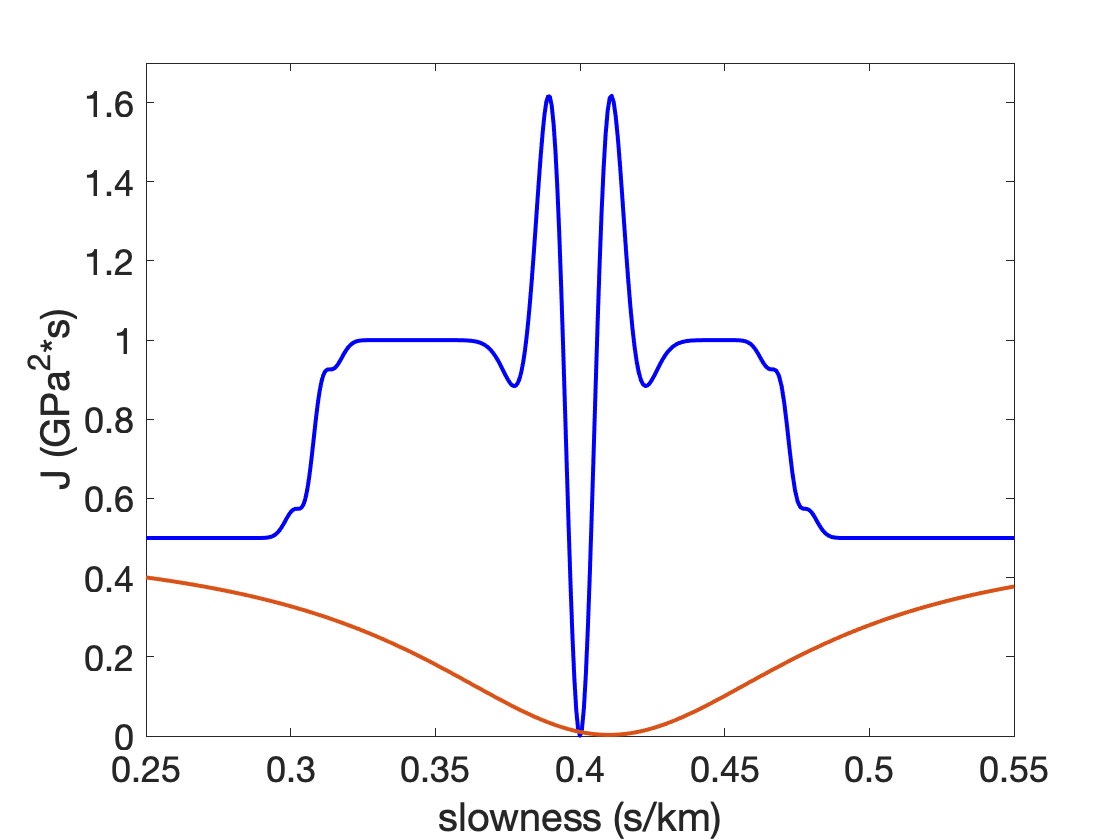}
\caption{Experiment 2a: reduced FWI (blue curve) and reduced ESI (red curve) objective functions versus slowness.}
\label{fig:asymmetric_fwi_vpm_01_40hz}
\end{figure}

\begin{figure}
\centering
\includegraphics[width=0.5\textwidth]{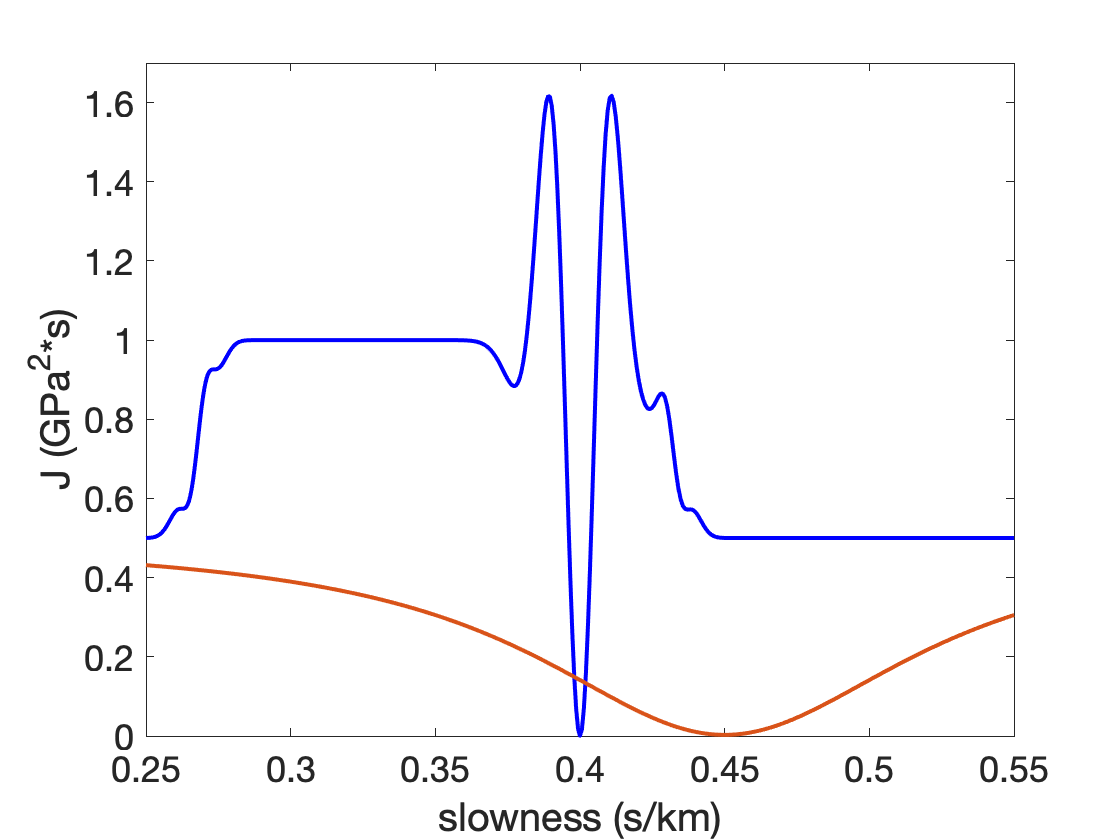}
\caption{Experiment 2b: reduced FWI (blue curve) and reduced ESI (red curve) objective functions versus slowness.}
\label{fig:asymmetric_fwi_vpm_05_40hz}
\end{figure}

\begin{figure}
\centering
\includegraphics[width=0.5\textwidth]{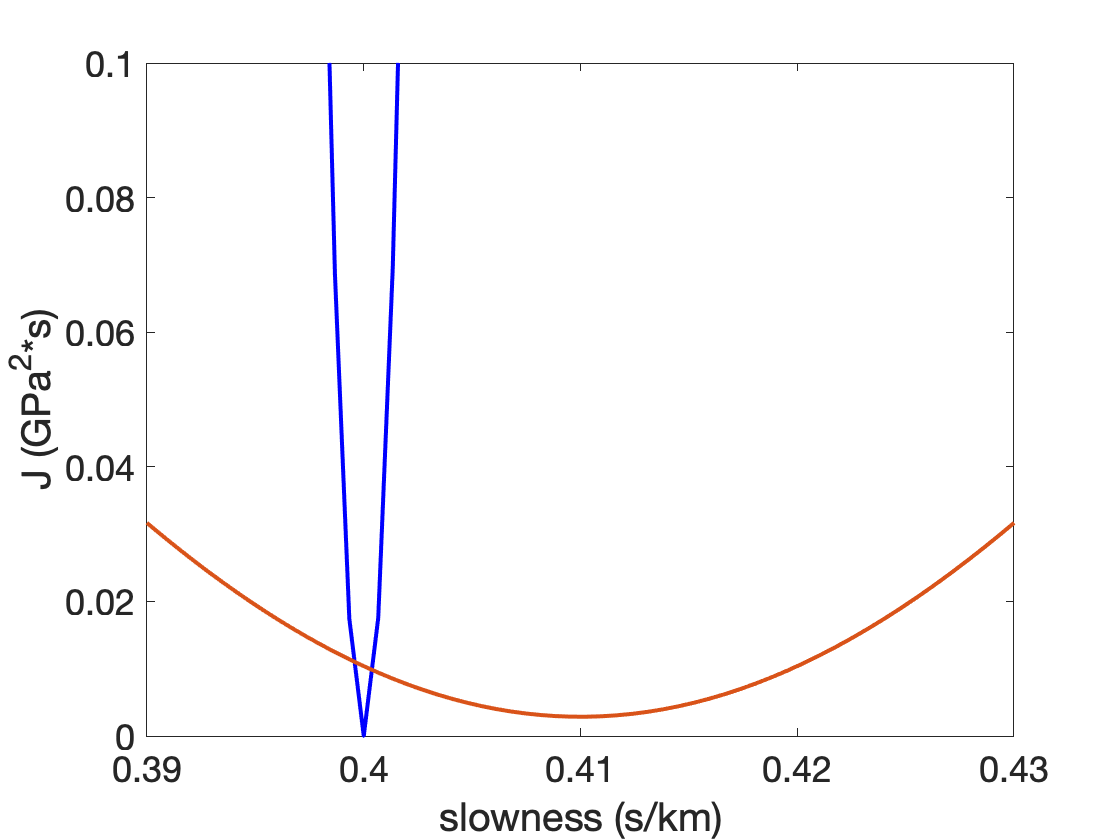}
\caption{A zoomed-in plot of Figure~\ref{fig:asymmetric_fwi_vpm_01_40hz} near the location of the global minimum.}
\label{fig:zoom_in_asymmetric_fwi_vpm_40hz}
\end{figure}

The target wavelet for experiment 2a is the wavelet of experiment 1 (Figure \ref{fig:figure_9}), shifted by 0.01 s, see Figure~\ref{fig:diff_target_wavelet}. The target wavelet for experiment 2b is shifted further, by 0.05 s. Both wavelets are supported in the interval $[-0.1, 0.1]$ s, so we use $\lambda = 0.1$ s here. The target slowness is still $m_*=0.4 $ s/km. The noise-free data for experiment 2a is shown in shown in Figure~\ref{fig:data_1}. The red curves in Figures~\ref{fig:asymmetric_fwi_vpm_01_40hz} and \ref{fig:asymmetric_fwi_vpm_05_40hz} show that in both cases, the estimated slowness $m$ is in error by less than $0.1$ s/km, consistent with Result \ref{thm:result2}. The error in the second case ($0.05$ s/km) is larger than the error in the first case (0.01 s/km, see Figure \ref{fig:zoom_in_asymmetric_fwi_vpm_40hz}). By changing the wavelet shape and piling up its energy near the support limit $\lambda=0.1$, we could obtain an error arbitrarily close to the bound \ref{eqn:merr}. That is, the bound in Result \ref{thm:result2} is sharp.

The blue curves in Figures \ref{fig:asymmetric_fwi_vpm_01_40hz} and \ref{fig:asymmetric_fwi_vpm_05_40hz} show the reduced FWI objective $\tilde{e}$. Definition \ref{thm:fwired} shows that the assumed maximum lag $\lambda$ plays a central role in the computation of $\tilde{e}$. The minimum is taken over all wavelets with support in $[-\lambda,\lambda]$. Result \ref{thm:result1} predicts that $\tilde{e} = 1/2$ for $|m-m_*| > 2 \lambda/r$ = 0.2 s/km for the choices made in these examples, and reference to the figures shows that this prediction holds.

\begin{figure}
\centering
\includegraphics[width=0.5\textwidth]{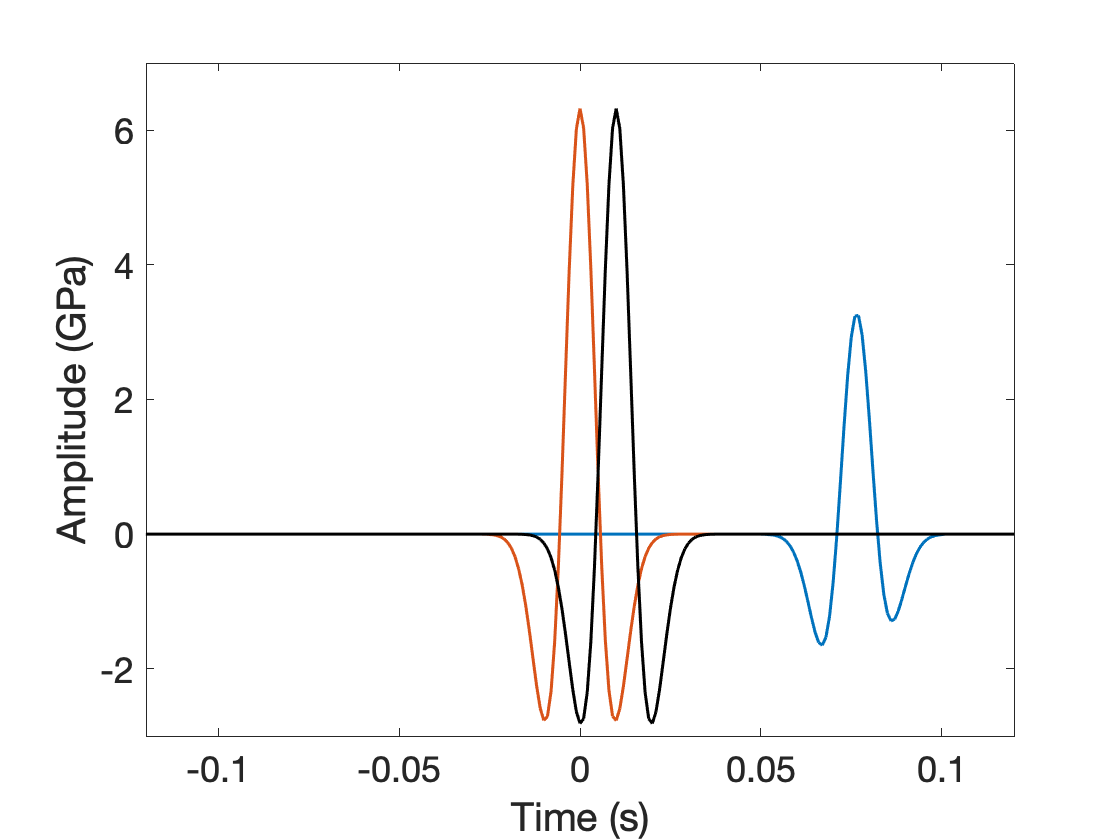}
\caption{Estimated wavelets for Experiment 2a: Blue curve is the initial wavelet at $m = 1/3$, $\alpha = 1$. Red curve is the estimated wavelet at $m = 0.41 $, $\alpha = 1$. The black curve is the target wavelet.}
\label{fig:ew_asy_fix_alpha}
\end{figure}

\begin{figure}
\centering
\includegraphics[width=0.5\textwidth]{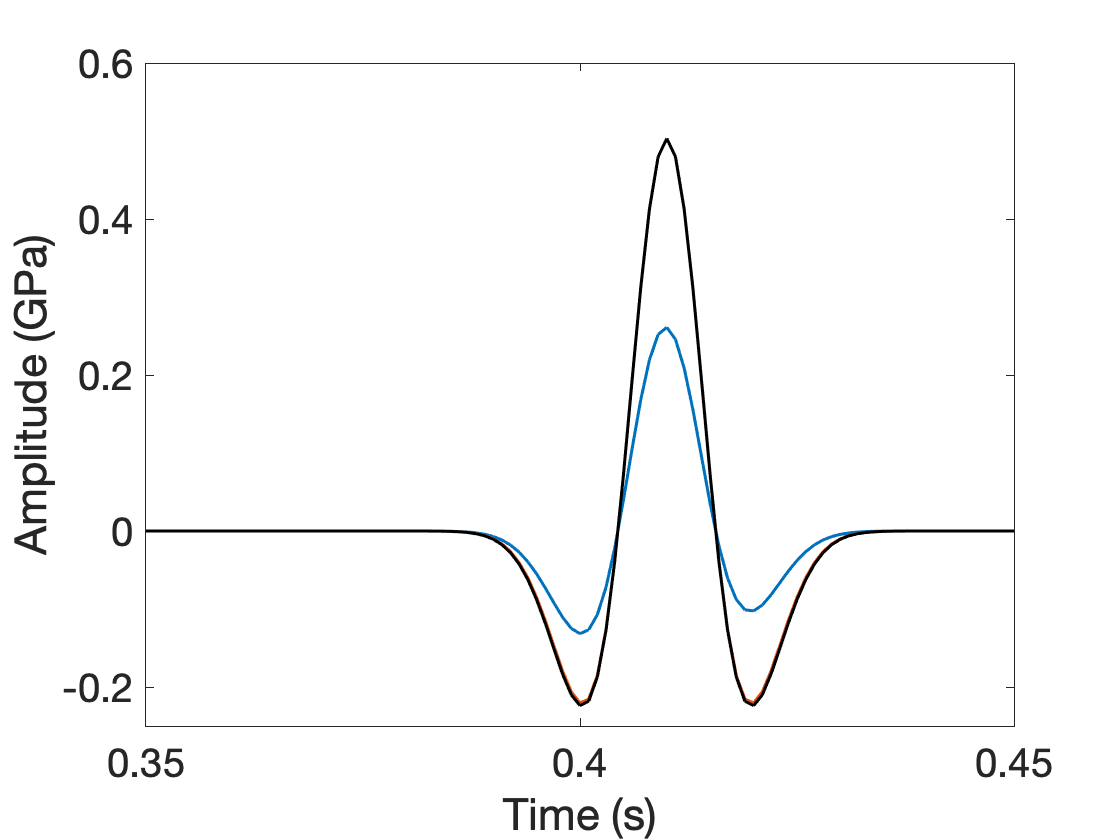}
\caption{Estimated data for Experiment 2a: Blue curve is the initial data at $m = 1/3$, $\alpha = 1$. Red curve is the estimated data at $m = 0.41 $, $\alpha = 1$. The black curve is the true data. Note that the red curve and black curve lay on top of eachother and are, therefore, difficult to distinguish.}
\label{fig:ed_asy_fix_alpha}
\end{figure}

\begin{figure}
\centering
\includegraphics[width=0.5\textwidth]{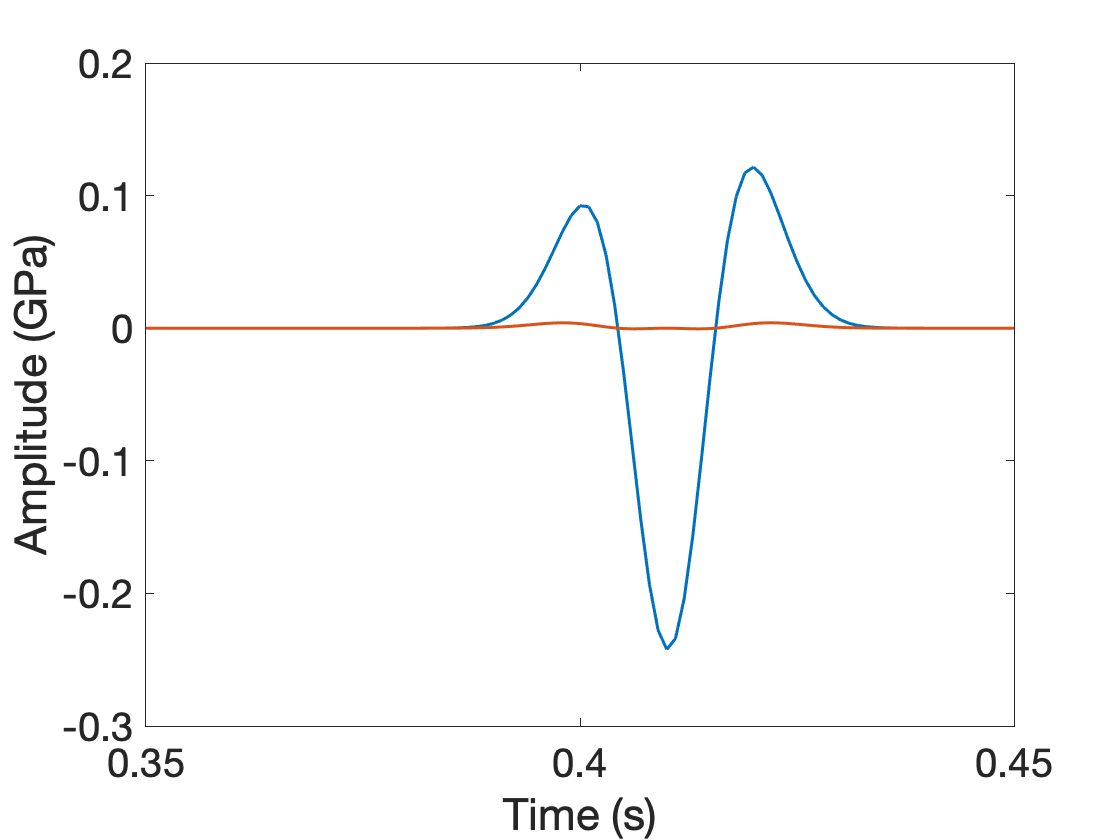}
\caption{Residual ($F[m]w[m;d] - d$) for Experiment 2a: Blue curve is the difference between the initial data and the true data. Red curve is the difference between the estimated data and the true data.}
\label{fig:residual_asy_fix_alpha}
\end{figure}

Minimization of the reduced ESI objective for the data of experiment 2a converged quickly to $m \approx 0.41$ from an initial estimate of $m=1/3$, and produced the wavelets shown in Figure \ref{fig:ew_asy_fix_alpha}. The predicted data also converged quickly to a final iterate closely matching the target data, as shown in Figures  \ref{fig:ed_asy_fix_alpha} and \ref{fig:residual_asy_fix_alpha}.

\begin{figure}
\centering
\includegraphics[width=.5\textwidth]{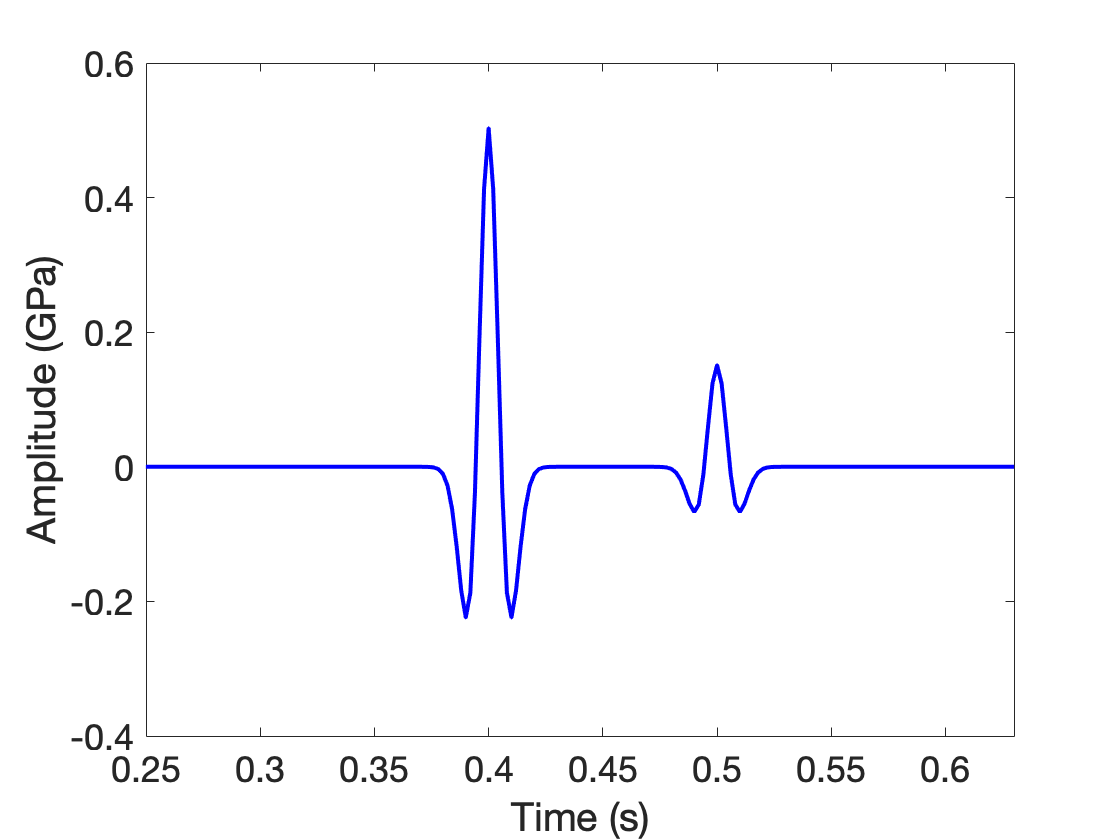}
\caption{Data for experiment 3 produced by adding to the data of Figure \ref{fig:nf_data} a shifted, scaled multiple of itself.} 
%(by 0.3) of a shifted (by 0.1 s) copy.}
\label{fig:noisydata_1}
\end{figure}

\noindent
\textbf{Experiment 3: 30\% Coherent Noise} Figure~\ref{fig:noisydata_1} displays the result of adding coherent noise to the data used in Experiment 1. This coherent noise is a small-amplitude shifted $40$ Hz Ricker wavelet centered at $t = 0.5$ s. The noise-to-signal ratio is $\eta = 0.3$. With this choice of $\eta$ and $\lambda=0.025$ we obtain an upper bound on the error from inequality~\eqref{eqn:merr}. Specifically,

$$|m-m_*| \leq \left(1+ \frac{2\eta(1+\eta)}{1-\eta(1+ \eta)}\right)\frac{\lambda}{r}$$
$$ \approx 0.057.$$ 

\begin{figure}
\centering
\includegraphics[width=.5\textwidth]{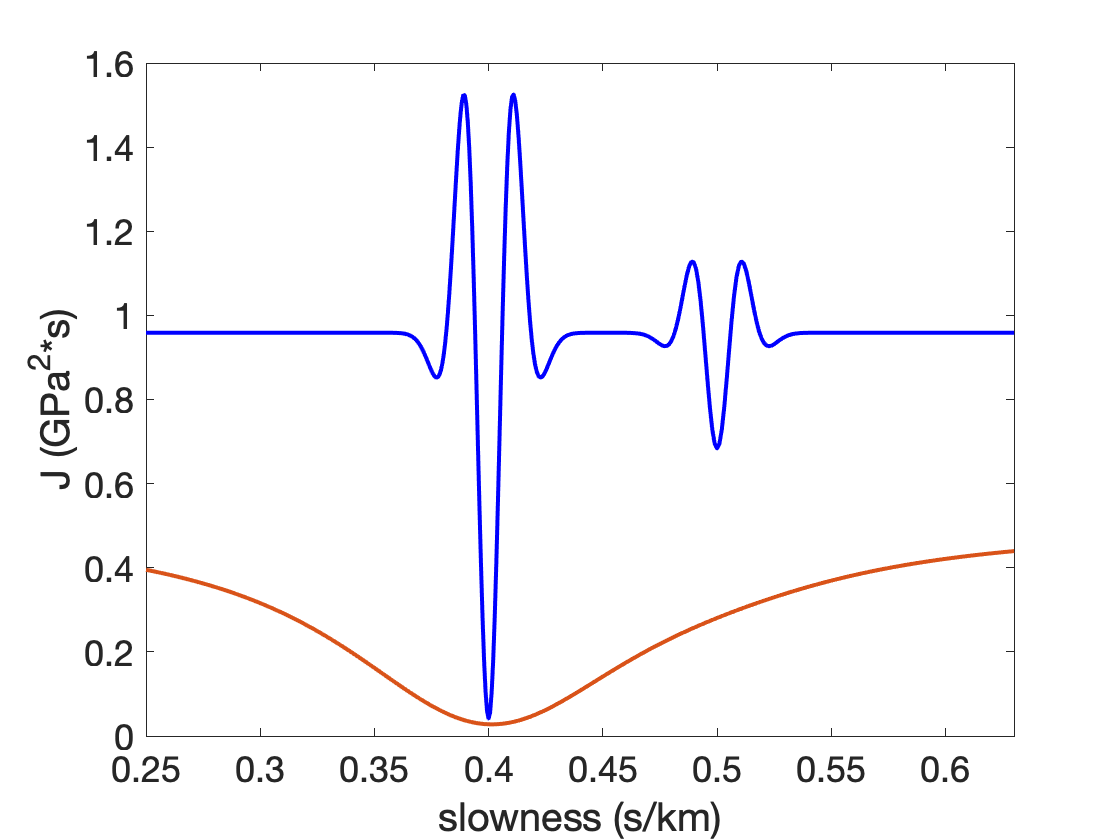}
\caption{Experiment 3: The reduced FWI and ESI objective functions versus slowness.}
\label{fig:figure_2}
\end{figure}

\begin{figure}
\centering
\includegraphics[width=.5\textwidth]{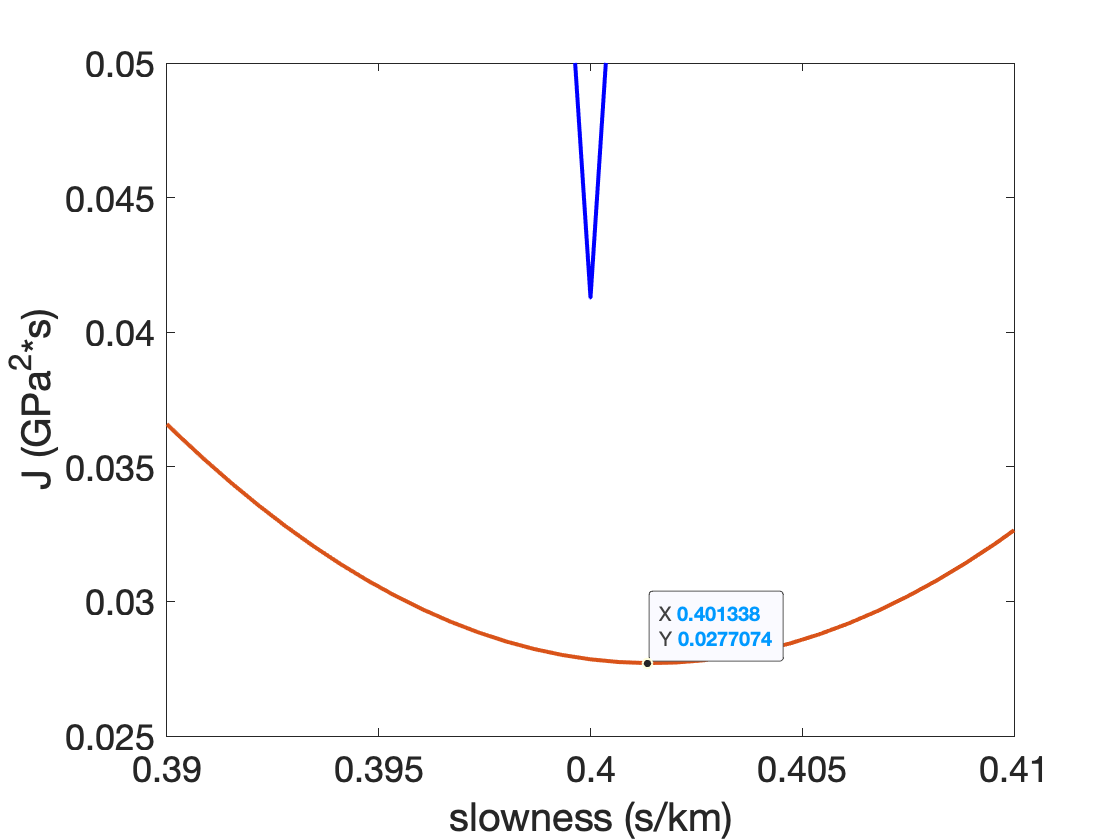}
\caption{Experiment 3: A zoomed-in plot of Figure~\ref{fig:figure_2} near the global minimum.}
\label{fig:figure_10}
\end{figure}

The minimum $m$ of the reduced ESI objective function is $\approx 0.401338$ s/km as shown in Figure~\ref{fig:figure_2} and seen more clearly in Figure~\ref{fig:figure_10}. The difference $|m-m_*| \approx 0.01338$ is substantially smaller than the upper error bound of $0.057$ predicted by Result \ref{thm:result2}.

\begin{figure}
\centering
\includegraphics[width=.5\textwidth]{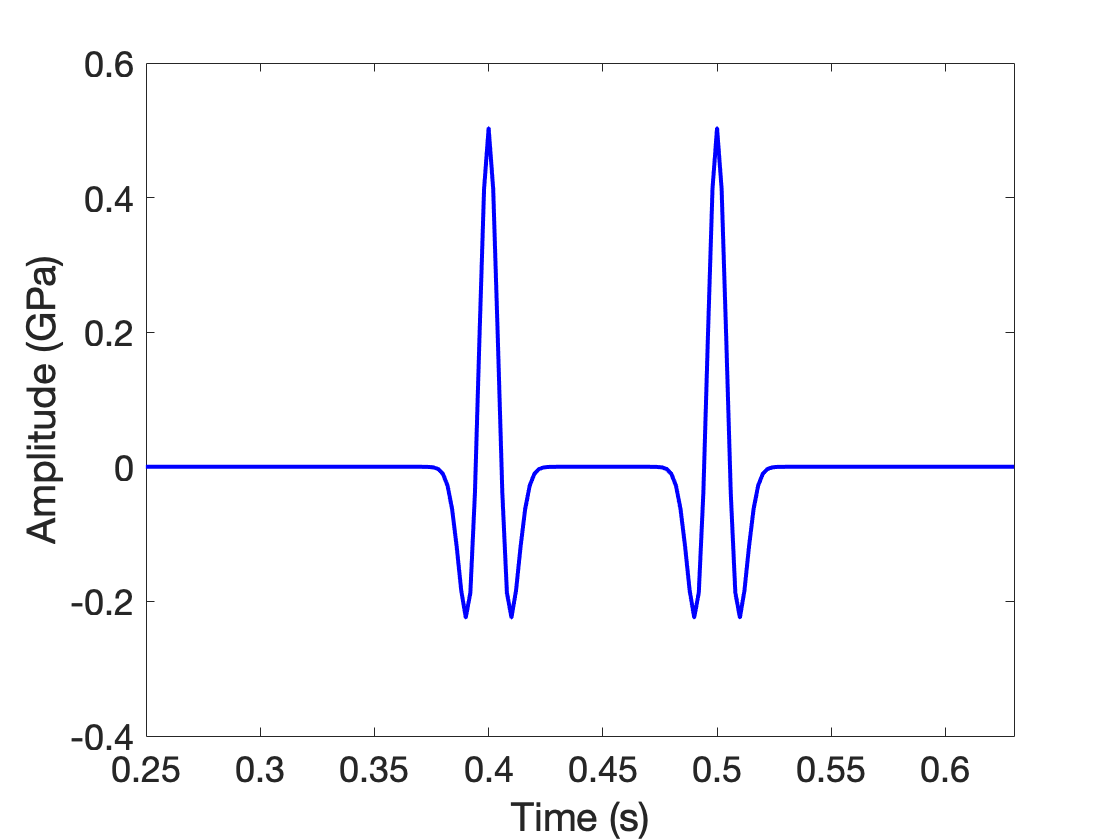}
\caption{Experiment 4: data. Signal is same as shown in Figure \ref{fig:nf_data} with the noise being a shifted copy of the signal.}
\label{fig:noisydata_2}
\end{figure}

\begin{figure}
\centering
\includegraphics[width=.5\textwidth]{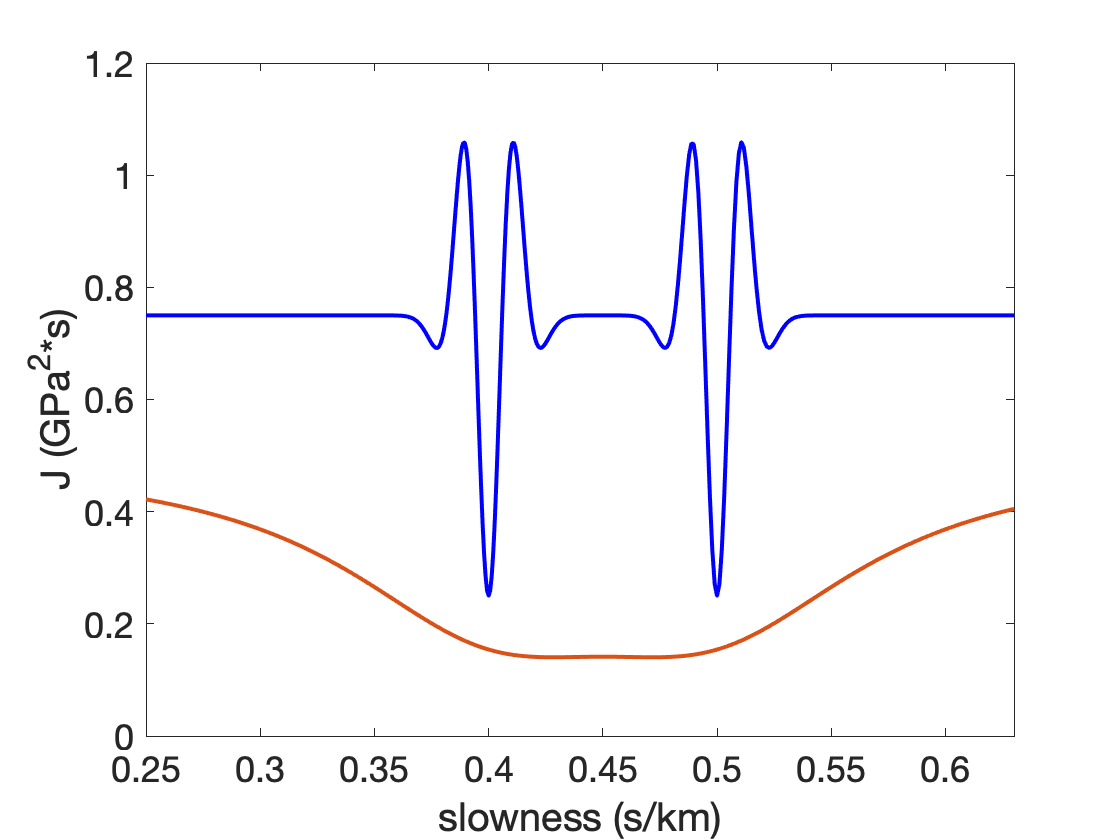}
\caption{Experiment 4: The restricted FWI (blue curve) and reduced ESI (red curve) objective functions versus slowness.}
\label{fig:figure_3}
\end{figure}

\noindent
\textbf{Experiment 4: 100\% Coherent Noise} \\
We have already pointed out in the previous section that no result like inequality~\eqref{eqn:merr} in Result \ref{thm:result2} could possibly hold for 100\% (or more) noise, since that condition includes the case of noise that annihilates the signal, making the reduced ESI objective constant in $m$. 
Still it seems worth illustrating this point in a different way. The next example uses data similar to that of Experiment 3 but with larger noise that violates the hypotheses of Result \ref{thm:result2}. The data is given in Figure~\ref{fig:noisydata_2}. The noise (centered at $t=0.5$ s) is a shifted copy of the signal. Thus the signal to noise ratio %($\eta$) 
is 1, which does not satisfy the assumption of Result \ref{thm:result2} that $\eta < \frac{\sqrt{5} - 1}{2}$. For this data, the reduced ESI objective function (with $\alpha=1$) shown in Figure~\ref{fig:figure_3} exhibits multiple local (and global) minima. Both the target slowness ($m_* = 0.4$ s/km) and the slowness $m = 0.5$ s/km, that would be inferred if the noise were the entire signal, are near local minimizers. 

Evidently this example could be modified, with similar results, by changing the time shift between the signal and the noise. Therefore the noise-to-signal ratio $\eta$ does not bound the difference between the target slowness $m_*$ and other stationary points of $\tJa$. That is, no bound like that stated in Result \ref{thm:result2} can possibly hold for noise levels this large.

\begin{figure}
\centering
\includegraphics[width=.5\textwidth]{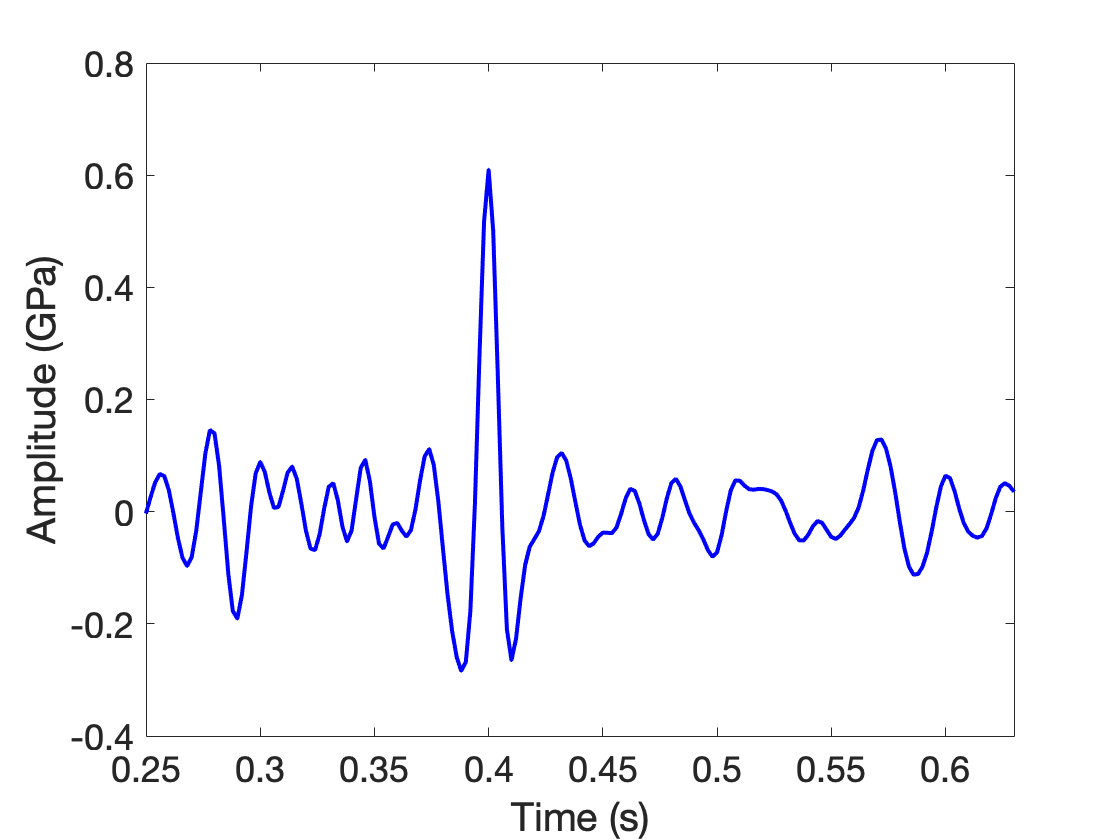}
\caption{Experiment 5: data.}
\label{fig:noisydata_3}
\end{figure}

\begin{figure}
\centering
\includegraphics[width=.5\textwidth]{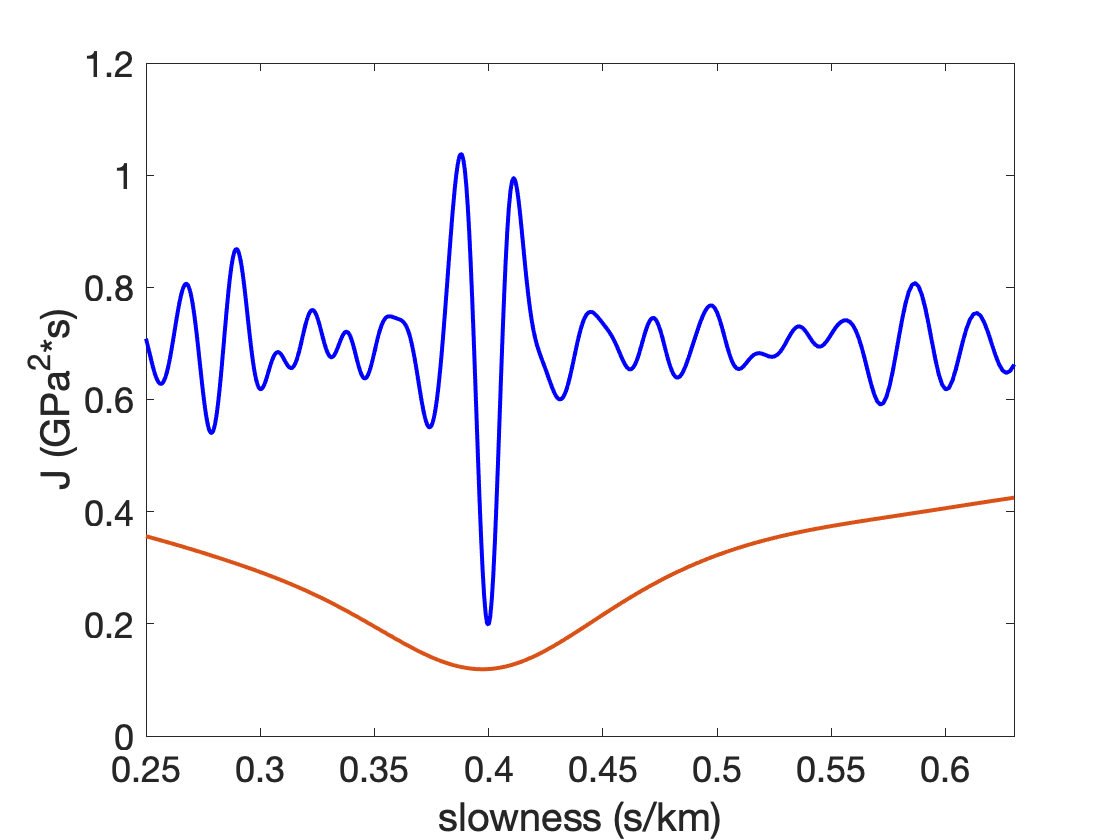}
\caption{Experiment 5: Restricted FWI (blue curve) and reduced ESI (red curve) objective functions plotted versus slowness.}
\label{fig:figure_4}
\end{figure}

%%%%%%%%%%%
%%%%%%%%%%%
\noindent
\textbf{Experiment 5: 100\% Filtered Random Noise} \\
We also mentioned in the last section that constraints on the type of noise may permit ESI to accurately estimate slowness and wavelet for noise-to-signal ratios higher than that permitted by Result \ref{thm:result2}.
Figure \ref{fig:noisydata_3} shows a data trace generated by adding to the signal in Figure \ref{fig:nf_data} a band-limited uniformly distributed random noise trace with $\eta = 1$, which again violates the approximately 60\% bound stated in Result \ref{thm:result2}. Inspection of the reduced ESI objective function (red curve in Figure \ref{fig:figure_4}) shows a unique local, and global, minimizer at almost precisely the target slowness used to generate the signal (Figure \ref{fig:nf_data}). The restricted FWI objective function (blue curve in Figure \ref{fig:figure_4}), on the other hand, has almost uniformly distributed local minima.

Despite their simplicity, these five experiments fully illustrate Result \ref{thm:result2}. More precisely, if the target wavelet $w_*$ vanishes for $|t|>\lambda$ and
the noise level $\eta$ is below $\approx 60\%$, the reduced ESI objective has all of its stationary points within $O(\lambda) + O(\eta)$ of the {\em FWI} global minimum.

\subsection{Discrepancy Algorithm and Solution of the Inverse Problem}

None of the examples so far have taken the final step of truncation to the prescribed maximum lag $\lambda$ necessary to produce a solution of the Inverse Problem as stated at the beginning of the last section. Also all of the examples discussed so far have employed a fixed, arbitrarily chosen penalty weight $\alpha=1$. Since the penalty operator (multiplication by $t$) penalizes spread of energy away from $t=0$, one would expect larger penalty weight to result in a smaller energy spread in the inverted wavelet, so less modification of the wavelet to conform to the maximum lag requirement of the Inverse Problem. The discrepancy algorithm increases $\alpha$ while maintaining data fit, so together with truncation of the wavelet at the prescribed lag (Result \ref{thm:result3}), should provide a route to solution of the inverse problem. The last two examples in this section illustrate 
Results \ref{thm:result3} and \ref{thm:result4}. 

\begin{figure}
\centering
\includegraphics[width=0.5\textwidth]{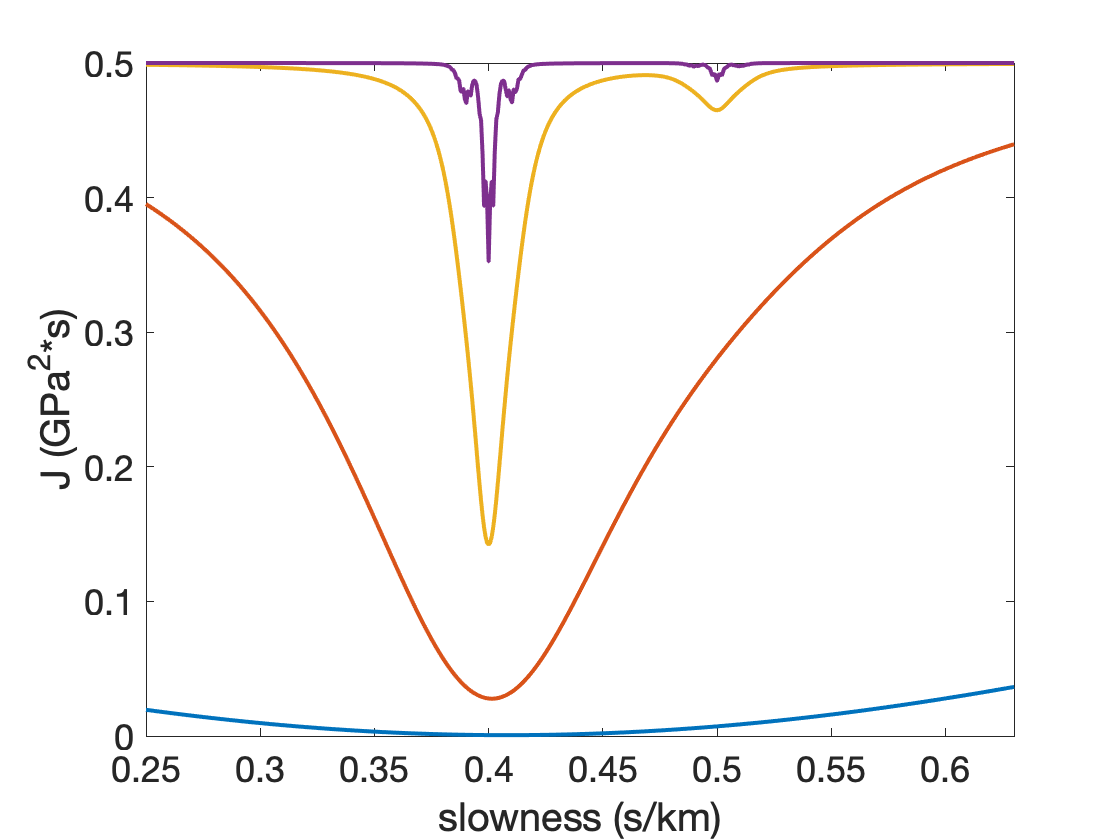}
\caption{The reduced ESI objective function corresponding to the data shown in Figure \ref{fig:noisydata_1} for various values of $\alpha$. The blue curve corresponds to $\alpha=0.1$. The red, yellow and purple curves
correspond to $\alpha = 1.0, 10.0$ and 100.0 respectively.}
\label{fig:noise_diff_alpha}
\end{figure}

We first discuss the selection of $\alpha$.
Figure~\ref{fig:noise_diff_alpha} shows the reduced ESI objective function with various $\alpha$ values for the data given in Figure~\ref{fig:noisydata_1}. The shape of the reduced ESI objective function changes drastically as $\alpha$ increases. For small $\alpha$ (blue curve in Figure~\ref{fig:noise_diff_alpha}), the objective function is very flat, but convex, with a unique stationary point. Achieving a small gradient (the usual stopping criterion for local optimization) for this choice of $\alpha$ will allow considerably more ambiguity in the estimated stationary point than does the same stopping criterion for a larger $\alpha$ (red, yellow or purple curves in Figure~\ref{fig:noise_diff_alpha}). This conclusion is partly captured in Result \ref{thm:result4}. Note that the yellow and purple curves exhibit local minimizers far from the target slowness, similar to reduced FWI. Evidently, $\alpha$ too small will inadequately constrain the slowness estimate, whereas $\alpha$ too large risks duplicating the behaviour of FWI.
However, note that the local minimum values for large $\alpha$, indeed all values of the objective, are quite large.
%for large $\alpha$. 
Because the penalty term $g$ decreases with $\alpha$, whereas the data error $e$ increases, the data error is the dominant contributor to these values for large $\alpha$.

The discrepancy algorithm adjusts $\alpha$ dynamically to improve the rate of convergence and the accuracy of the model estimates obtained by ESI, while maintaining a prescribed level of data error. A proper constraint on data error prevents $\alpha$ from being increased beyond the level where multiple local minima appear. 
%The examples to come apparently use such a proper constraint; as will be elaborated in the Discussion section, 
%We are not currently aware of an algorithm for choosing such a constraint, nor of a theoretical basis for such a choice. 

\noindent
The discrepancy algorithm requires that an acceptable range of data error ($[e_-, e_+]$) be prescribed.
We use a relatively large range, based on a target signal-to-noise of $3$. The actual signal-to-noise ratio in this example is $\approx 3.33$, corresponding to $e \approx 0.045$. We choose the error bounds for the discrepancy algorithm as $e_- \approx 0.027, e_+ \approx 0.11$. 
Brent's method \cite[]{Brent:TCJ71}, used for the $m$ update in the discrepancy algorithm, requires these inputs:
\begin{itemize}
    \item[1. ]Initial search interval $[m_{\rm min}, m_{\rm max}]$: We use $m_{\rm min} = 0.33, m_{\rm max}=0.65$ for every $m$ update cycle;
    \item[2. ]A stopping tolerance for the gradient length: We use (absolute) tolerances of $0.01$ for Experiment 6, $0.001$ for Experiment 7.
\end{itemize}

%Brent's method, as described in \cite[]{Brent:TCJ71}, also uses a separate stopping criterion, based on a minimum steplength. However we disabled this branch in our implementation of the algorithm.

\noindent
\textbf{Experiment 6: Coherent Noise}\\
In this example, we use the data depicted in Figure~\ref{fig:noisydata_1}.
%together with the discrepancy algorithm. Brent's method for updating $m$ (by finding a root of $\nabla_m \tJa$) stops when the gradient length drops below 1\% of its initial value.
The progress of the algorithm is described in Tables~\ref{table:coherent_alpha_first}-\ref{table:coherent_m_sec}. Table~\ref{table:coherent_alpha_first} shows the update of $\alpha$ using the rule \ref{eqn:alphasecant} from the initial value $\alpha=0$, attaining a value of $e$ in the prescribed range after three iterations. Table~\ref{table:coherent_m_first} shows the first Brent's method update of $m$ from its initial value of $0.343$ s/km. It ends with an approximate stationary point of $\tJa$ having been found, and indeed a very good estimate of $m_*$. However the data error $e$ has dropped below $e_-$, so another $\alpha$ update is necessary. As shown in Table ~\ref{table:coherent_alpha_sec}, a single application of rule \ref{eqn:alphasecant} is sufficient to bring $e$ into compliance with the required bounds. Another 14 iterations of Brent's method (Table ~\ref{table:coherent_m_sec}) compute another approximate stationary point of $\tJa$. This time, the data error $e$ also remains in the prescribed range, so the algorithm terminates. 

\begin{figure}
\centering
\includegraphics[width=.5\textwidth]{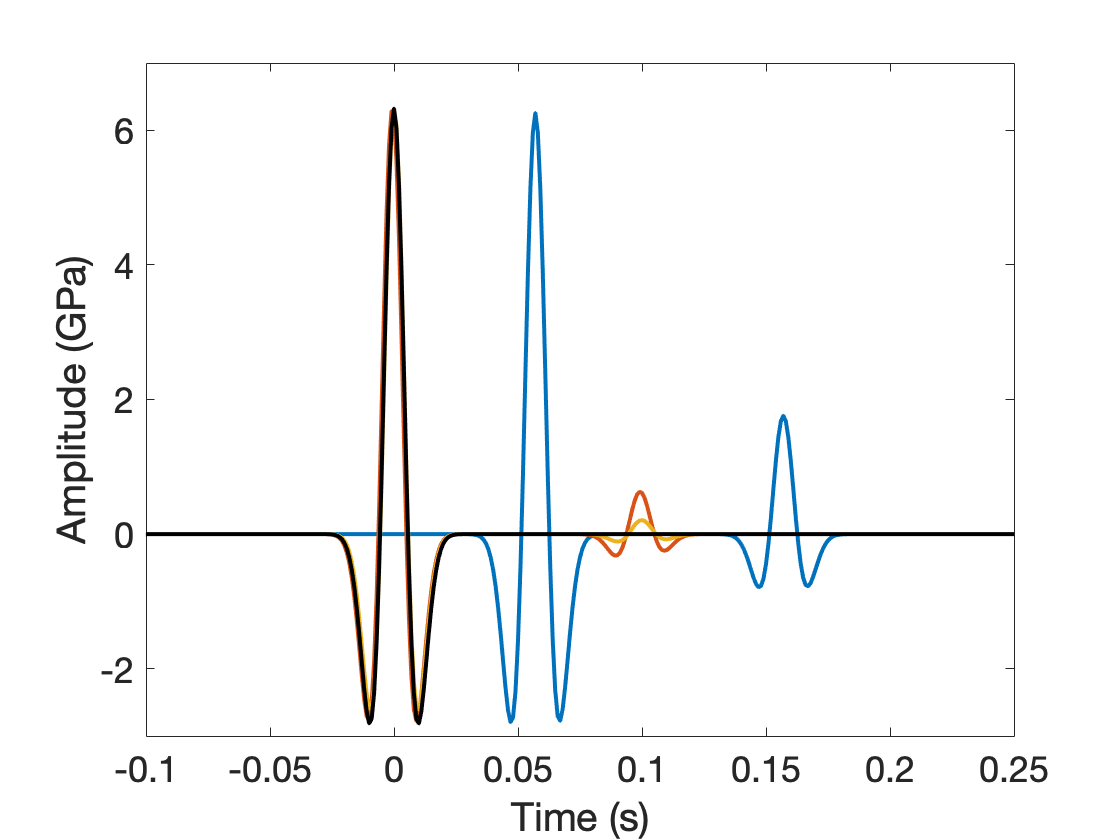}
\caption{Experiment 6: Estimated wavelets $w_\alpha[m;d]$. Blue curve is the initial wavelet at $m = 0.343$, $\alpha =0.284184$. Red curve is the estimated wavelet after the first update of $m$: $m = 0.401035$, $\alpha = 1.136737$. Yellow curve is the estimated wavelet after the second (final) update of $m$: $m = 0.400113$, $\alpha = 2.273473$. The black curve is the target wavelet.}
\label{fig:ew5}
\end{figure}

\begin{figure}
\centering
\includegraphics[width=.5\textwidth]{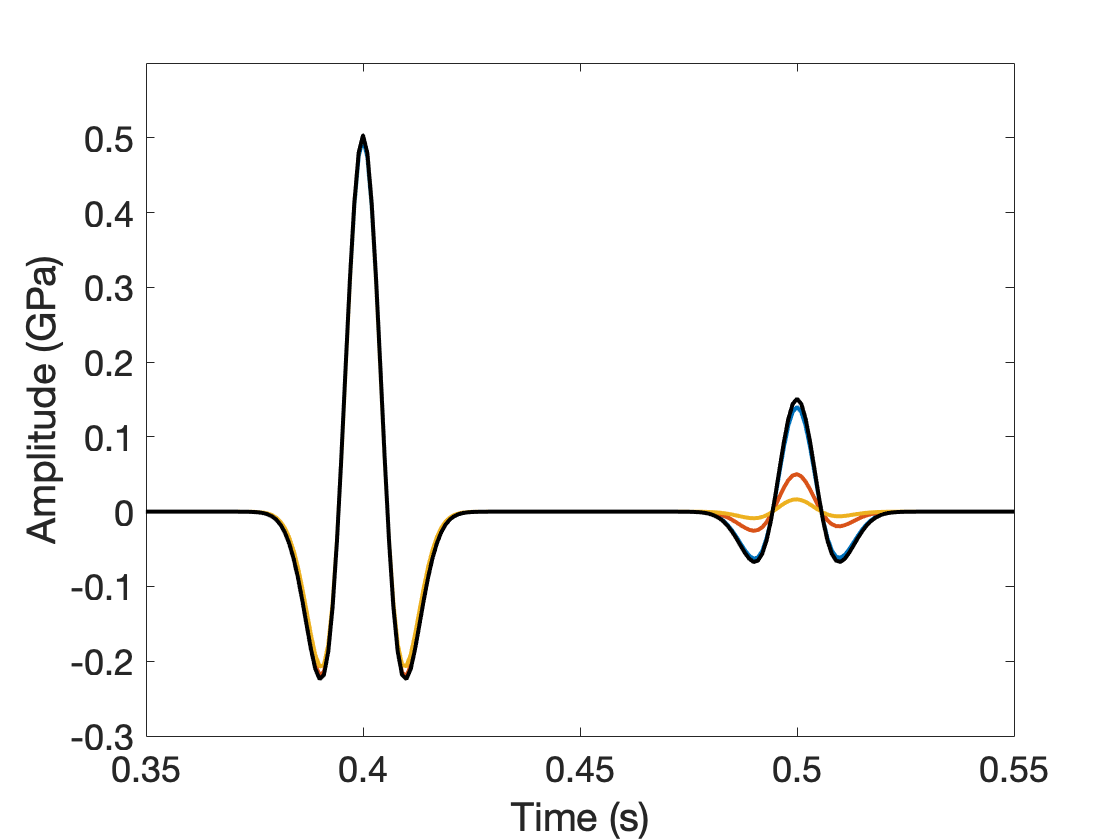}
\caption{Experiment 6: estimated data $F[m]w_\alpha[m;d]$. Blue curve is the initial data at $m = 0.343$, $\alpha =0.284184$.  Red curve is the estimated data after the first update of $m$: $m = 0.401035$, $\alpha = 1.136737$. Yellow curve is the estimated data after the second (final) update of $m$: $m = 0.400113$, $\alpha = 2.273473$. The black curve is the true data.}
\label{fig:ed5}
\end{figure}

\begin{figure}
\centering
\includegraphics[width=.5\textwidth]{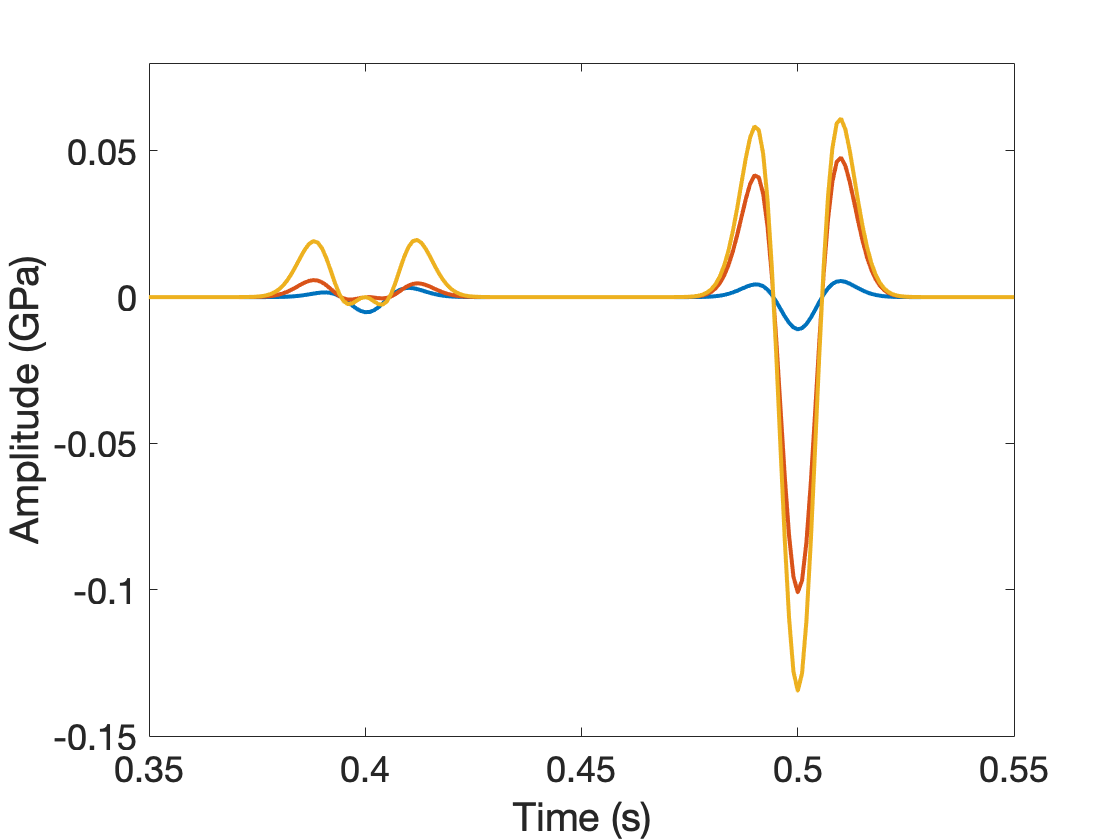}
\caption{Experiment 6: Data Residual ($F[m]w_\alpha[m;d]- d$). Blue curve is the difference between the initial data at $m = 0.343$, $\alpha =0.284184$ and the true data. Red curve is the difference between the estimated data after the first update of $m$: $m = 0.401035$, $\alpha = 1.136737$ and the true data. Yellow curve is the difference between the estimated data after the second (final) update of $m$: $m = 0.400113$, $\alpha = 2.273473$ and the true data.}
\label{fig:res5}
\end{figure}

Recall that the estimated wavelet $w_\alpha[m,d](t)$ (see Figure~\ref{fig:ew5}) is a solution of the normal equation given by Equation \eqref{eqn:norm}. The estimated data shown in Figure~\ref{fig:ed5} is computed by $F[m]w_\alpha$ and its residual is shown in Figure~\ref{fig:res5}.

\begin{comment}
\begin{table*}
\begin{minipage}{80mm}
  \begin{center}
\vspace{0.3in}
    \begin{tabular}{c|c}
    \hline\hline
      \textbf{Parameters: }  & \\
      \hline 
            Frequency & $40$ Hz\\
            \hline 
            target slowness & 0.4 s/km\\
             \hline
             target snr ($snr_{tgt}$)  & 3\\
              \hline
              true snr & $1/0.3 \approx3.33$\\
              \hline
              $\|d\|$ & 0.045445\\  \hline
               $e_{tgt} = 0.5*(1/snr_{tgt})^2$ & 0.055556  \\  \hline
    %           $C$ & 0.7\\ \hline
                $\gamma_{-}$ &$ (0.7)^2 = 0.490000$ \\  \hline
                  $\gamma_{+}$ & $(0.7)^{-2} = 2.040816 $\\   \hline
                 $e_{-} = \gamma_{-} e_{tgt}$ &  0.027222 \\  \hline
                  $e_{+}= \gamma_{+} e_{tgt}$ &  0.113379 \\  \hline
                  gradient tolerance & $10^{-2}$\\  \hline
                    max number of iterations allowed for the outer loop & 100\\      \hline    
                 initial slowness: m & 0.343 s/km\\     \hline
  %                  2 starting points for the Brent's method& $a = 1/3$, $b = 0.65$\\     
      \hline
    \end{tabular}
    \caption{Experiment 4, parameters used for the discrepancy algorithm.}
    \label{table:coherent_parameters}
  \end{center}
  \end{minipage}
\end{table*}
\end{comment}

\begin{table*}
\begin{minipage}{75mm}
  \begin{center}   
    \begin{tabular}{c|c|c|c|c}
    \hline\hline
      \textbf{iteration: } & {$\alpha$} & {$g$}& {$e$} \\
           \hline 
  1  &   0.284184  &.  0.371103&      0.003140 \\
  2  &   0.568368  &   0.311447 &     0.022460  \\
  3  &   1.136737  &   0.204342 &     0.102216  \\         \hline
    \end{tabular}
    \caption{Experiment 6, discrepancy algorithm: $\alpha$ updates for initial $m = 0.343$. Initial $\alpha = 0$.}
     \label{table:coherent_alpha_first}
  \end{center}
  \end{minipage}
\end{table*}

\begin{comment}
\begin{table*}
\begin{minipage}{90mm}
  \begin{center}
\vspace{0.3in}
    \begin{tabular}{c|c|c|c|c|c|c}
    \hline\hline
      \textbf{$i$: }  &   \textbf{$g$}  & \textbf{$e$}  & \textbf{$m$}& \textbf{$J_{\alpha}$}& \textbf{$\nabla J_{\alpha}$}  \\ 
      \hline 
   1   &  0.035018 &   0.403247    &   0.622695 &   0.448496   &  0.463686\\
  2    &  0.089906  &  0.140974      &  0.478014   & 0.257147  &   2.614541 \\
  3    &  0.011959  &  0.017344    &   0.405674  &  0.032797  &   0.803269\\
  4    &  0.009643 &   0.018478  &    0.400642 &   0.030938  & -0.070100\\
    \hline
    \end{tabular}
    \caption{Experiment 6, discrepancy algorithm: updates of $m$ with after first update of $\alpha = 1.136737$. Initial $m=0.343$.}
     \label{table:coherent_m_first}
  \end{center}
  \end{minipage}
\end{table*}
\end{comment}

\begin{table*}
\begin{minipage}{90mm}
  \begin{center}
\vspace{0.3in}
    \begin{tabular}{c|c|c|c|c|c|c}
    \hline\hline
      \textbf{$i$: }  &   \textbf{$g$}  & \textbf{$e$}  & \textbf{$m$}& \textbf{$\tJa$}& \textbf{$\frac{d}{dm} \tJa$}  \\ 
      \hline 
  1   &   0.035018   &   0.403247    &  0.622695   &   0.448496     &  0.463686\\
  2   &   0.089906  &    0.140974   &   0.478014   &   0.257147    &   2.614541\\
  3   &   0.011959   &    0.017344  &    0.405674  &    0.032797   &    0.803269\\
  4   &   0.028659   &   0.025577   &   0.381536  &    0.062608   &    -3.049986\\
  5   &   0.009643   &   0.018478   &   0.400642  &    0.030938   &    -0.070100\\
  6   &   0.010393   &   0.017888   &   0.403158  &    0.031317   &    0.370812\\
  7   &   0.009914   &   0.018178   &   0.401900   &    0.030989  &     0.151012\\
  8   &   0.009752   &   0.018327   &   0.401271   &   0.030929   &    0.040569\\
  9   &   0.009691   &   0.018402   &   0.400956  &    0.030925   &    -0.014743\\
 10   &   0.009720   &   0.018364  &    0.401114  &    0.030924   &    0.012919\\
 11   &   0.009705   &   0.018383  &    0.401035  &    0.030924    &   -0.000911\\
    \hline
    \end{tabular}
    \caption{Experiment 6, discrepancy algorithm: updates of $m$ after first update of $\alpha = 1.136737$. Initial $m=0.343$.}
     \label{table:coherent_m_first}
  \end{center}
  \end{minipage}
\end{table*}

\begin{table*}
\begin{minipage}{75mm}
  \begin{center}  
    \begin{tabular}{c|c|c|c|c}
    \hline\hline
     \textbf{iteration: } & {$\alpha$} & {$g$}& {$e$}& {$m$} \\
           \hline 
  1&     2.273473   & 0.009705      &0.033737    &  0.401035\\
     \hline
    \end{tabular}
    \caption{Experiment 6, discrepancy algorithm: second update of $\alpha$ after first update of $m=0.401035$. Initial $\alpha =1.136737 $.}
     \label{table:coherent_alpha_sec}
  \end{center}
  \end{minipage}
\end{table*}

\begin{table*}
\begin{minipage}{90mm}
  \begin{center}
\vspace{0.3in}
    \begin{tabular}{c|c|c|c|c|c|c}
    \hline\hline
       \textbf{$i$: }  & \textbf{$g$}  & \textbf{$e$}  & \textbf{$m$}& \textbf{$\tJa$}& \textbf{$\frac{d}{dm} \tJa$}  \\ 
             \hline 
   1&    0.002303  &   0.475832 &    0.637763 &    0.487735 &     0.114887\\
  2 &    0.011948  &   0.336897 &     0.485548 &    0.398651 &     0.700990\\
  3 &    0.007396  &   0.037167 &     0.409441 &    0.075396 &       5.288562\\
  4 &    0.020844  &   0.111823 &     0.371387 &    0.219561  &    -7.541345\\
  5 &    0.007280  &   0.040128 &     0.390414 &    0.077754 &     -5.521535\\
  6 &    0.002986  &   0.033854 &     0.399927 &    0.049290 &     -0.128092\\
  7 &    0.004197 &    0.034122 &     0.404684&     0.055816 &     2.827748\\
  8 &    0.003301 &    0.033728 &     0.402306&     0.050789 &     1.382890\\
  9 &    0.003068  &   0.033732 &     0.401116 &    0.049590 &     0.631507\\
 10 &    0.003008 &    0.033779 &     0.400522 &    0.049327 &      0.252196\\
 11 &    0.002992 &    0.033813 &     0.400225 &    0.049280 &     0.062106\\
 12 &    0.002988  &   0.033833 &     0.400076 &    0.049278 &     -0.032988\\
 13 &    0.002990  &   0.033823 &     0.400150 &    0.049278 &     0.014561\\
 14 &    0.002989  &   0.033828 &     0.400113 &    0.049278 &     -0.009213\\
        \hline
    \end{tabular}
    \caption{Experiment 6, discrepancy algorithm: second update of $m$ after second update of $\alpha=2.273473$. Initial $m= 0.401035$. At the last iteration, the gradient norm termination criterion for Brent's method is satisfied, and the data error $e$ is within the prescribed bounds, so the algorithm terminates.}
     \label{table:coherent_m_sec}
  \end{center}
  \end{minipage}
\end{table*}

Note that the estimated wavelet in Figure~\ref{fig:ew5} does not satisfy the maxiumum lag constraint posed in the Inverse Problem. As explained in the last section, the second step of our algorithm for solution of the inverse problem is the truncation of the wavelet to a specified maximum lag. Result \ref{thm:result3} gives estimates of both lag $\lambda$ and data error $\epsilon$, in terms of the maximum lag $\mu$ of the target wavelet $w_*$ and the noise-to-signal ratio $\eta$, for which the conditions of the Inverse Problem should be satisfied by the truncated wavelet $w[m;d]$ at a stationary point $m$ of $\tJa$. 

\begin{figure}
\centering
\includegraphics[width=.5\textwidth]{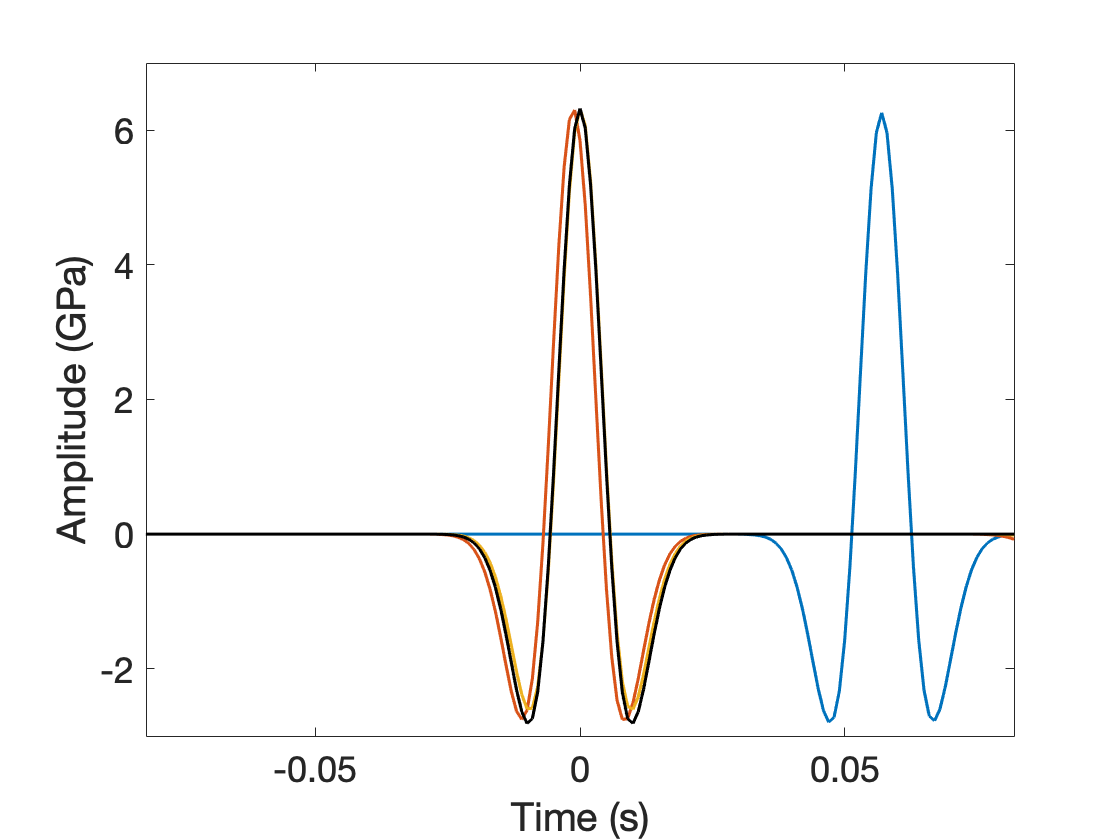}
\caption{Experiment 6: Truncated estimated wavelets ${\bf 1}_{[-\lambda,\lambda]}w_\alpha[m;d]$ for several iterates of the discrepancy algorithm. Blue curve is the initial wavelet at $m = 0.343$, $\alpha =0.284184$. Red curve is the estimated wavelet after the first update of $m$: $m = 0.401035$, $\alpha = 1.136737$. Yellow curve is the estimated wavelet after the second (final) update of $m$: $m = 0.400113$, $\alpha = 2.273473$. The black curve is the target wavelet.}
\label{fig:cutted_w_coh03}
\end{figure}

\begin{figure}
\centering
\includegraphics[width=.5\textwidth]{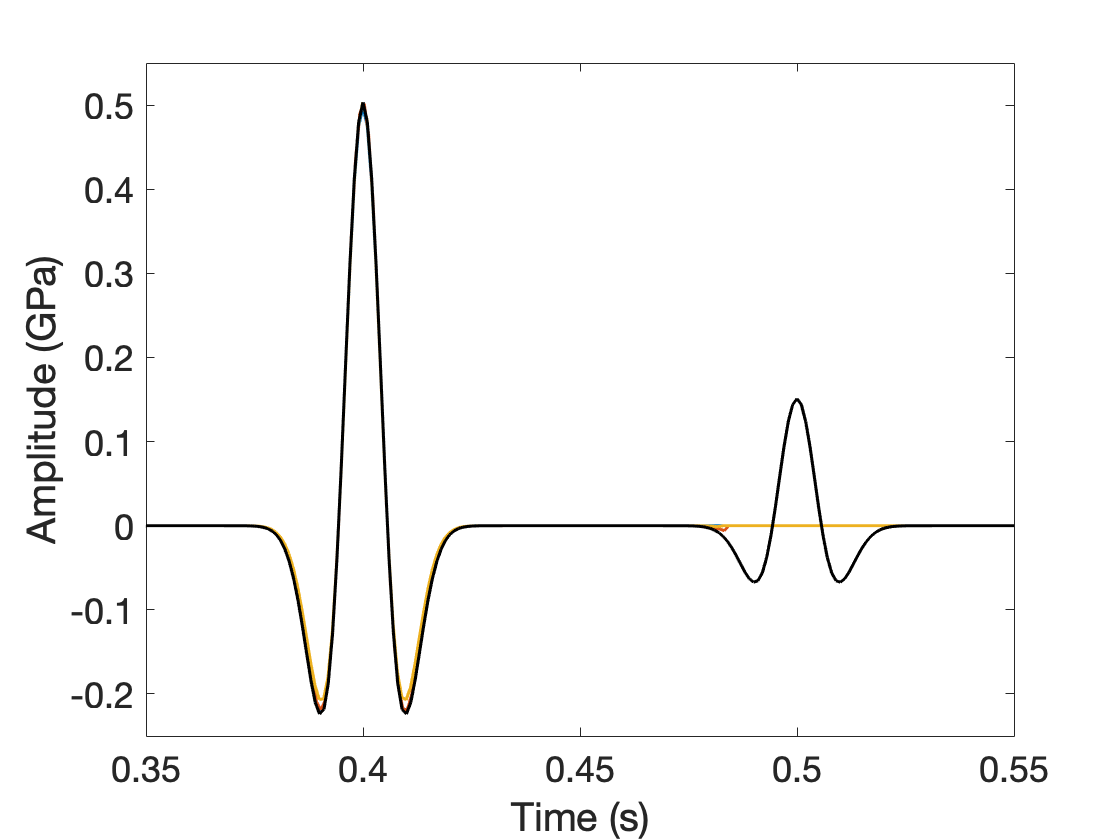}
\caption{Experiment 6: predicted data from truncated wavelets, $F[m]{\bf 1}_{[-\lambda,\lambda]}w_\alpha[m;d]$ for several iterates of the discrepancy algorithm. Blue curve is the initial data at $m = 0.343$, $\alpha =0.284184$.  Red curve is the estimated data after the first update of $m$: $m = 0.401035$, $\alpha = 1.136737$. Yellow curve is the estimated data after the second (final) update of $m$: $m = 0.400113$, $\alpha = 2.273473$. The black curve is the target data.}
\label{fig:cutted_data_coh03}
\end{figure}

\begin{figure}
\centering
\includegraphics[width=.5\textwidth]{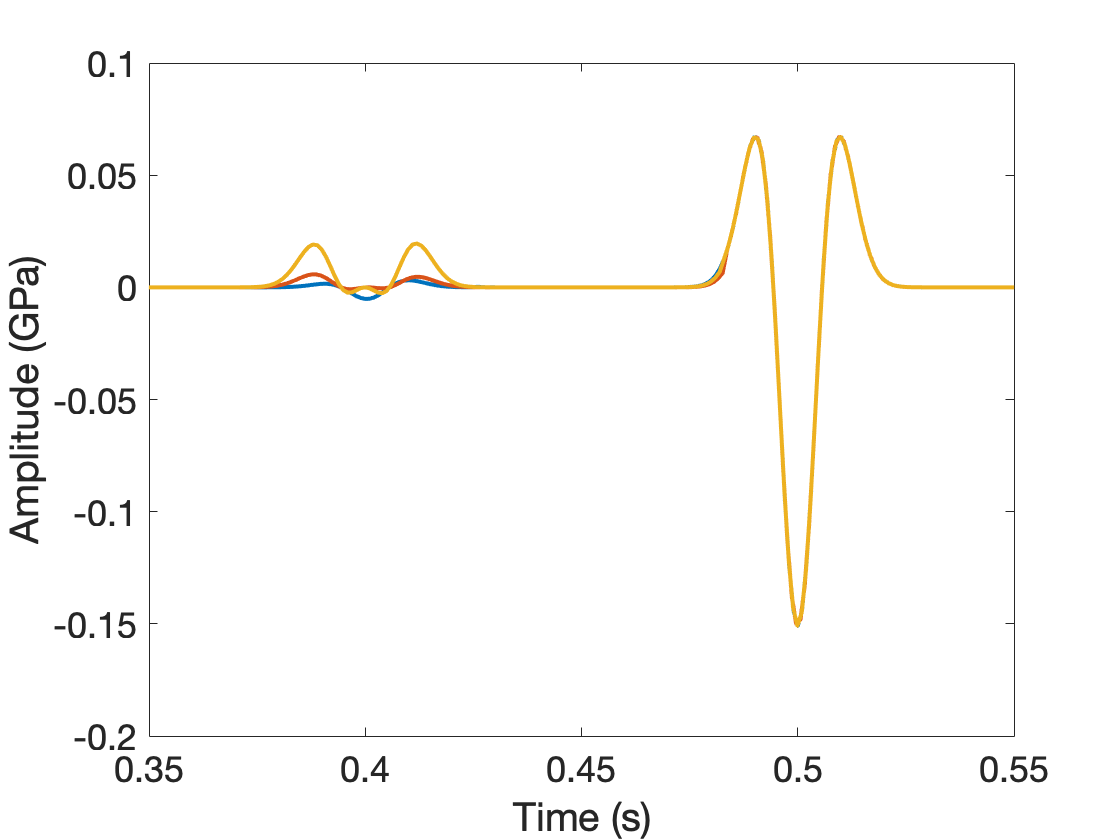}
\caption{Experiment 6: data residual ($F[m]{\bf 1}_{[-\lambda,\lambda]}w_\alpha - d$) for several iterates of the discrepancy algorithm. Blue curve is the residual at the initial values $m = 0.343$, $\alpha =0.284184$. Red curve is the residual after the first update of $m$: $m = 0.401035$, $\alpha = 1.136737$. Yellow curve is the residual after the second (final) update of $m$: $m = 0.400113$, $\alpha = 2.273473$.}
\label{fig:cutted_res_coh03}
\end{figure}

In this synthetic example, we can use the known values  $\mu=0.025$ and $\eta=0.3$ in these estimates. However both estimates result from repeated use of the triangle inequality, and are very conservative. The estimated lag $\lambda$ from inequality \eqref{eqn:lambbd} is $\lambda = 0.082$, which turns out to be usable. Then the estimated wavelet is truncated to $[-\lambda, \lambda] = [-0.082, 0.082]$. This turns out to be a reasonable estimate for the truncation interval, as can be seen in plots of the truncated estimated wavelet (Figure~\ref{fig:cutted_w_coh03}), the data obtained by pairing this wavelet with the estimated slowness (Figure~\ref{fig:cutted_data_coh03}), and the resulting data residual (Figure~\ref{fig:cutted_res_coh03}).

On the other hand, the estimate \eqref{eqn:epsbd} is quite pessimistic, estimating data error at nearly 100\%. Instead, we simply compute the data error obtained with the estimated slowness and the truncated wavelet just described. We obtain
\[
\epsilon = \frac{\|F[m]{\bf 1}_[-\lambda,\lambda]w-d\|}{\|d\|} \approx 0.29,
\]
which is only slightly greater than the value obtained by substituting $m=m_*, w=w_*$. 

Therefore we tentatively conclude that one should  estimate a reasonable maxiumum lag $\lambda$ via equation \eqref{eqn:lambbd}, truncate the wavelet, and then obtain the relative data error directly from the resulting predicted data, rather than relying on estimate \ref{eqn:epsbd}. 

\begin{figure}
\centering
\includegraphics[width=.5\textwidth]{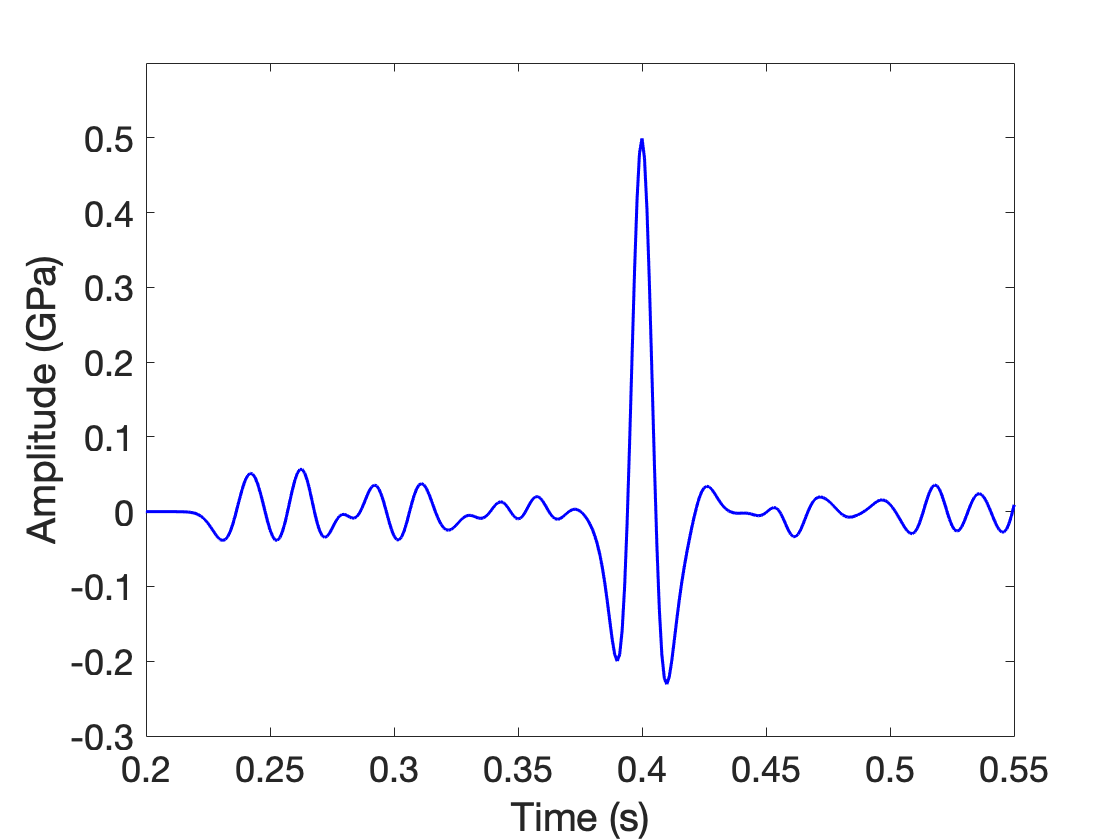}
\caption{Experiment 7: data.}
\label{fig:noisydata_4}
\end{figure}

%%%%%%%%%%%
%%%%%%%%%%%
\noindent
\textbf{Experiment 7: Filtered Random Noise} \\ The data for this experiment is shown in Figure~\ref{fig:noisydata_4}. It consists of the noise-free data shown in Figure \ref{fig:nf_data}, contaminated with uniformly distributed random noise, filtered by the noise-free wavelet and scaled to have noise-to-signal ratio $\eta=0.3$, as in the previous example. 

\begin{figure}
\centering
\includegraphics[width=.5\textwidth]{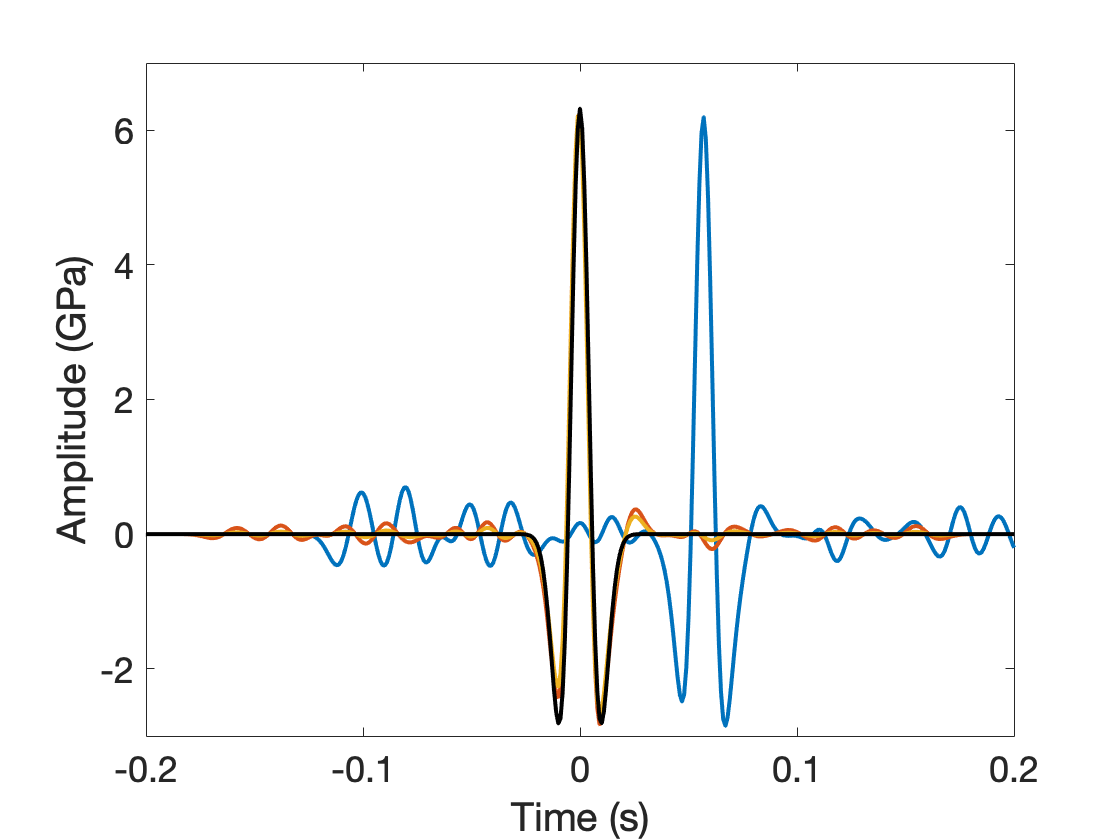}
\caption{Experiment 7: Estimated wavelets $w_\alpha[m;d]$ for several iterations of the discrepancy algorithm. Blue curve is the initial wavelet at $m = 0.343$, $\alpha =0.251370$. Red curve is the estimated wavelet after the first update of $m$: $m = 0.400444$, $\alpha = 1.005480$. Yellow curve is the estimated wavelet after the second (final) update of $m$: $m = 0.400499$, $\alpha = 2.010960 $. The black curve is the target wavelet.}
\label{fig:rn_snr333_ew}
\end{figure}

\begin{figure}
\centering
\includegraphics[width=.5\textwidth]{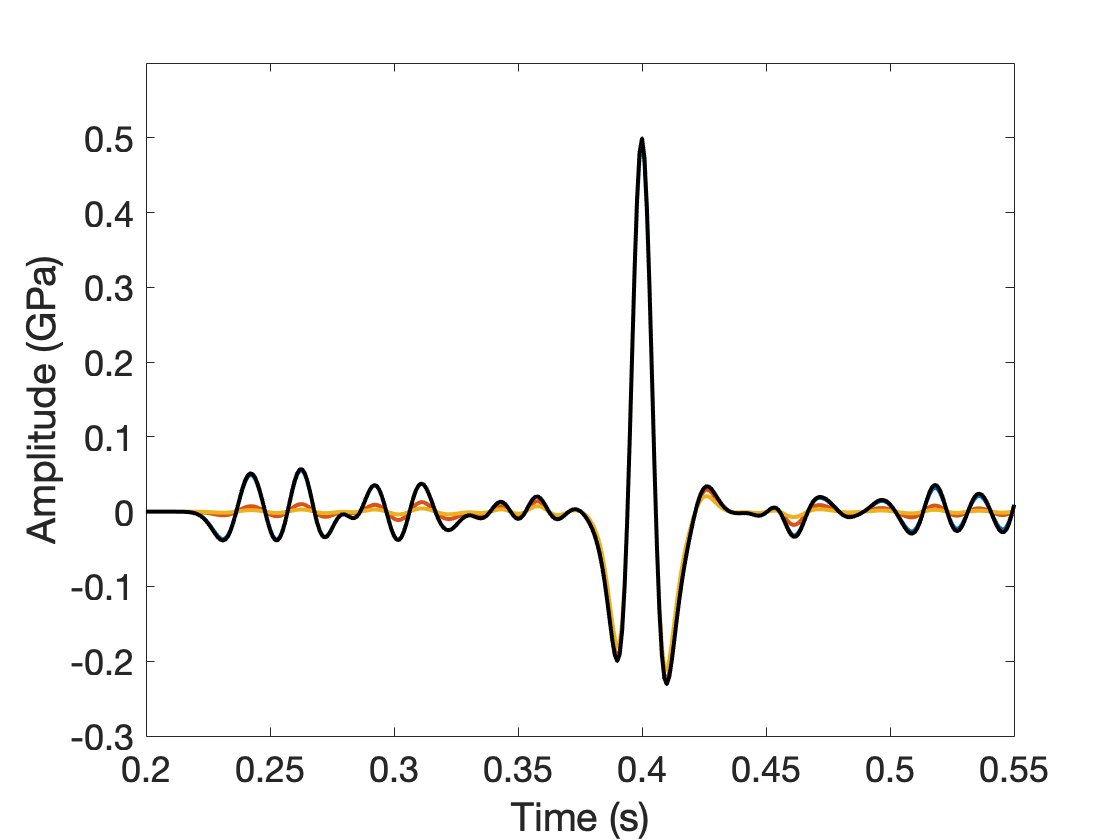}
\caption{Experiment 7: Estimated data $F[m]w_\alpha[m;d].$ Blue curve is the initial data at $m = 0.343$, $\alpha =0.251370$. Red curve is the estimated data after the first update of $m$: $m = 0.400444$, $\alpha = 1.005480$. Yellow curve is the estimated data after the second (final) update of $m$: $m = 0.400499$, $\alpha = 2.010960 $. The black curve is the target data.}
\label{fig:rn_snr333_ed}
\end{figure}

\begin{figure}
\centering
\includegraphics[width=.5\textwidth]{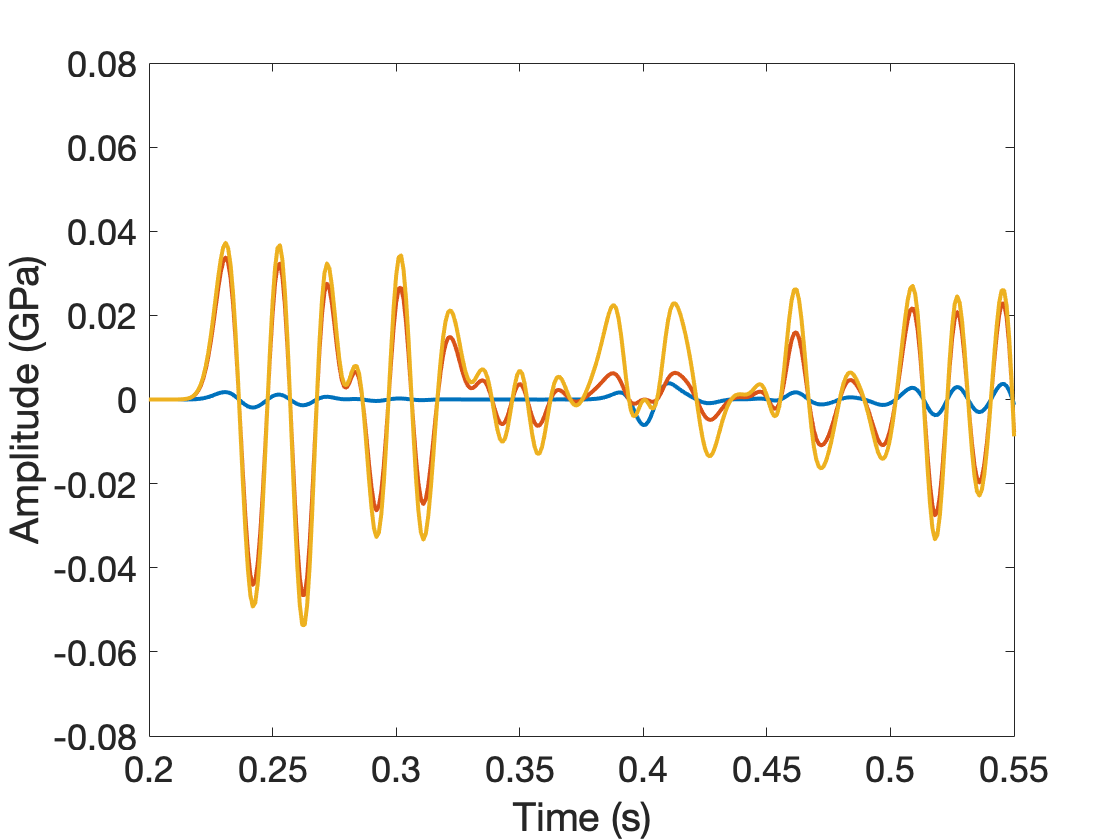}
\caption{Experiment 7: Data residual ($F[m]w_\alpha[m;d] - d$) for several iterations of the discrepancy algorithm. Blue curve is the residual at initial $m = 0.343$, $\alpha =0.251370$. Red curve is the residual after the first update of $m$: $m = 0.400444$, $\alpha = 1.005480$. Yellow curve is the residual after the second (final) update of $m$: $m = 0.400499$, $\alpha = 2.010960 $.}
\label{fig:rn_snr333_res}
\end{figure}

The course of the discrepancy algorithm is similar to that in the previous example, so we do not show the iterations in detail. Both $\alpha$ and $m$ are updated twice, at which point a stationary point of $\tJa$ has been found for which the bounds on $e$ are satisfied. The estimated wavelet extends quite far from $t=0$. The  
%Tables~\ref{table:rna_alpha_first}-\ref{table:rnd_m_sec} show iterations of the algorithm.
Figures~\ref{fig:rn_snr333_ew},~\ref{fig:rn_snr333_ed} and~\ref{fig:rn_snr333_res} show the wavelets, predicted data, and residuals respectively, computed in the course of two iterations. 
%Note that $\alpha$ is really a weight parameter here. That is, smaller $\alpha$ leads to heavier weight on the data term $d(t+mr)$ as we compute the wavelet. That's why our initial wavelet starting from small $\alpha$ fits the data so well (overfitting), leading to a small residual.

\begin{figure}
\centering
\includegraphics[width=.5\textwidth]{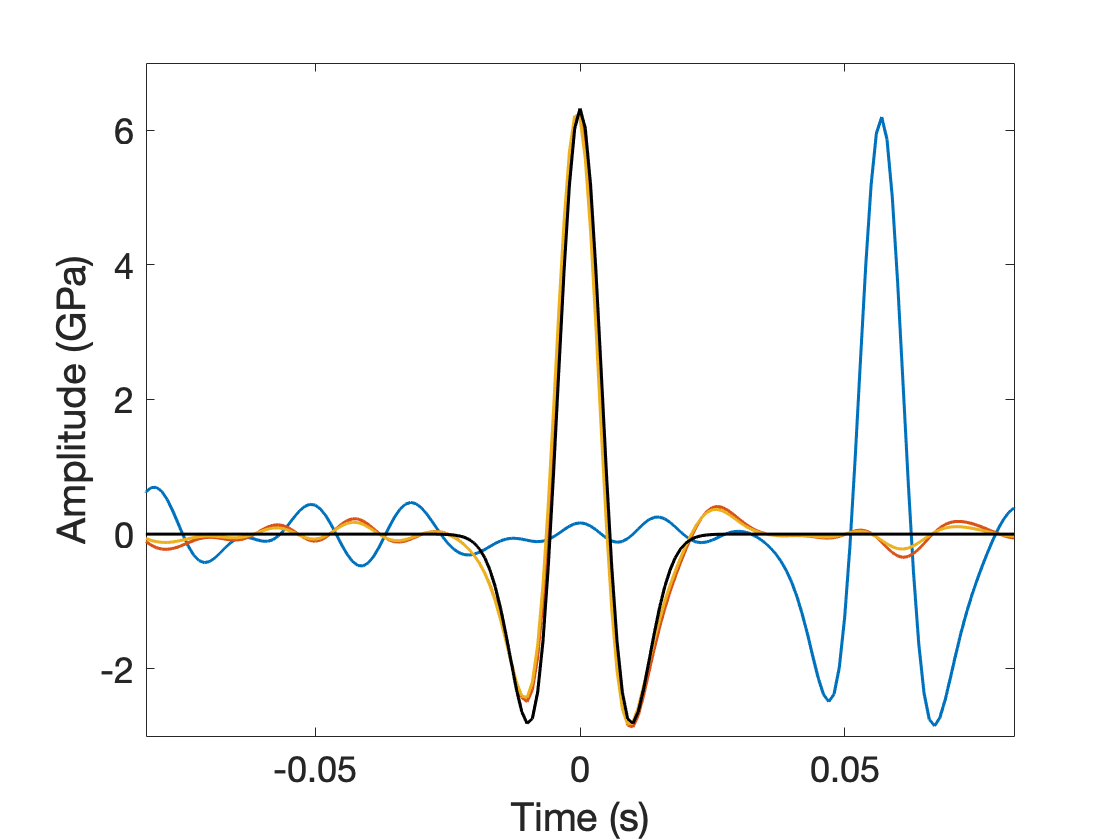}
\caption{Experiment 7: Truncated estimated wavelets ${\bf 1}_{[-\lambda,\lambda]}w_\alpha[m;d]$, $\lambda=0.082$ for several iterations of the discrepancy algorithm. Blue curve is truncated wavelet at initial $m = 0.343$, $\alpha =0.251370$. Red curve is the truncated wavelet after the first update of $m$: $m = 0.400444$, $\alpha = 1.005480$. Yellow curve is the truncated wavelet after the second (final) update of $m$: $m = 0.400499$, $\alpha = 2.010960 $. The black curve is the target wavelet.}
\label{fig:cutted_rn_snr333_ew}
\end{figure}

\begin{figure}
\centering
\includegraphics[width=.5\textwidth]{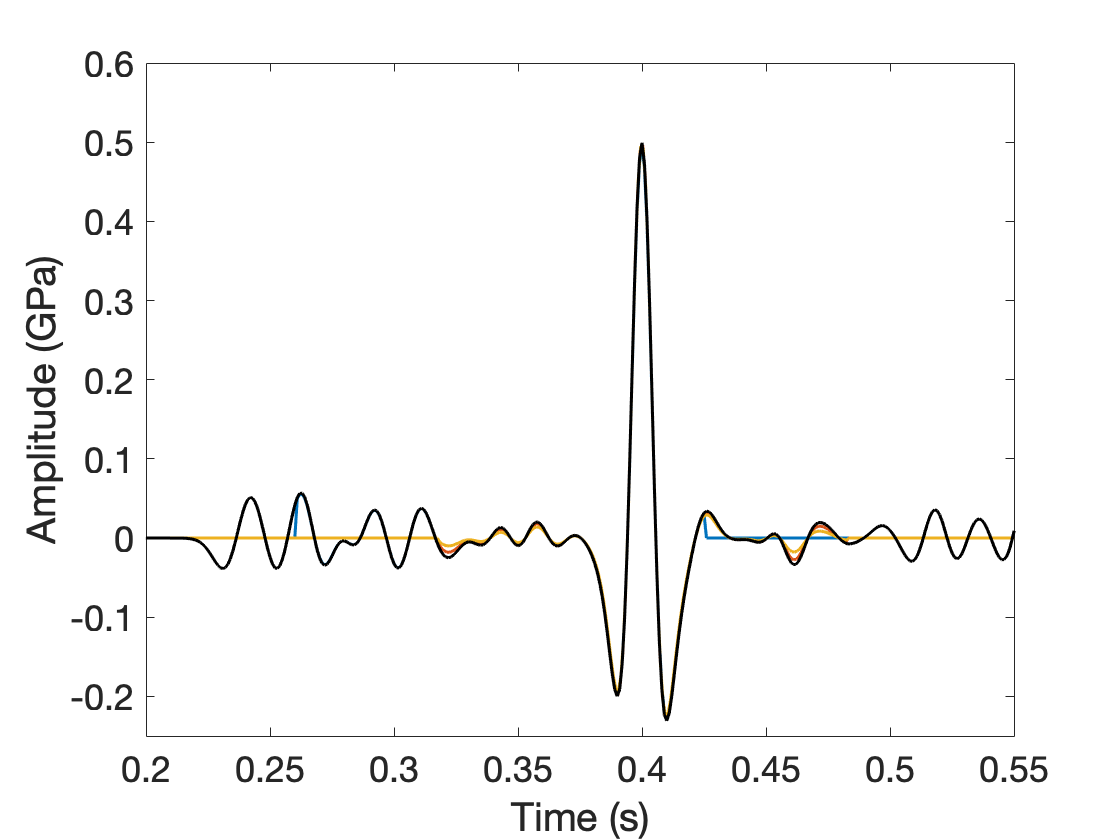}
\caption{Experiment 7: Estimated data from truncated wavelets, $F[m]{\bf 1}_{[-\lambda,\lambda]}w_\alpha[m;d]$ for several iterations of discrepancy algorithm. Blue curve is the initial data at $m = 0.343$, $\alpha =0.251370$. Red curve is the estimated data after the first update of $m$: $m = 0.400444 $, $\alpha = 1.005480$. Yellow curve is the estimated data after the second (final)  update of $m$: $m = 0.400499$, $\alpha = 2.010960 $. The black curve is the target data.}
\label{fig:cutted_rn_snr333_ed}
\end{figure}

\begin{figure}
\centering
\includegraphics[width=.5\textwidth]{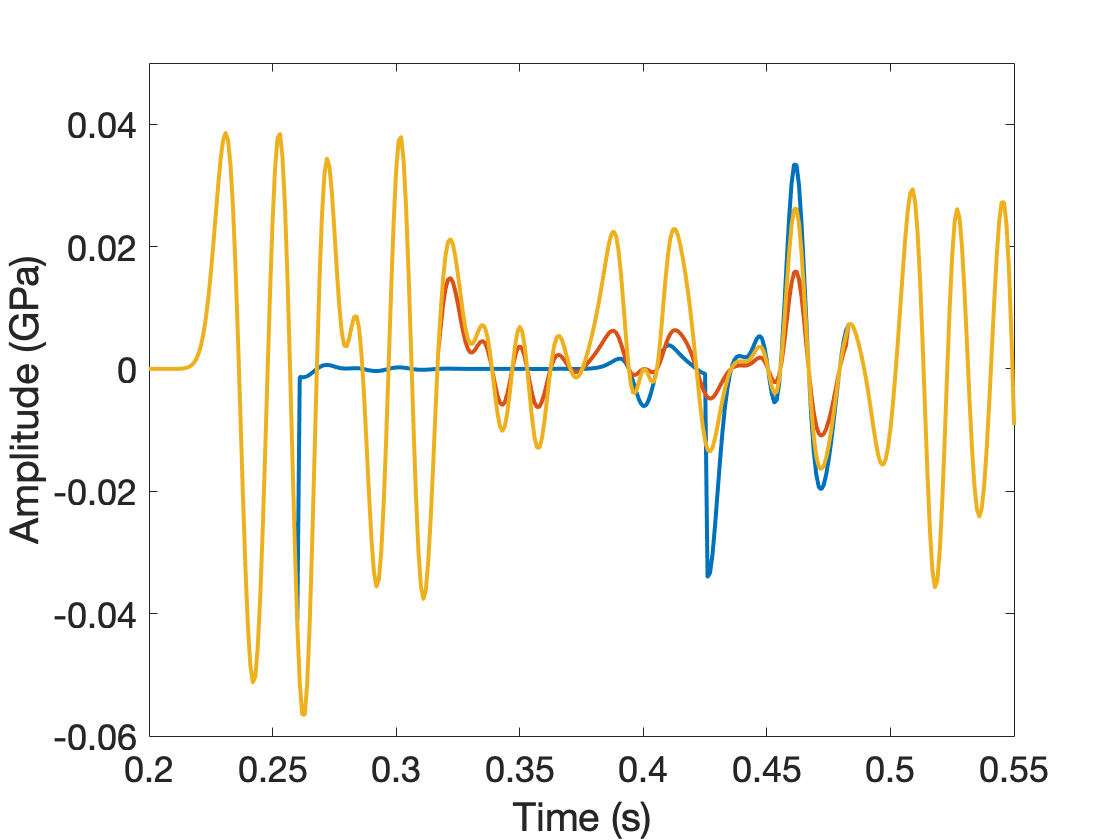}
\caption{Experiment 7: Data residual using truncated wavelets ($F[m]{\bf 1}_{[-\lambda,\lambda]}w_\alpha[m;d] - d$) for several iterations of the discrepancy algorithm. Blue curve is initial residual at $m = 0.343$, $\alpha =0.251370$. Red curve is residual after the first update of $m$: $m = 0.400444$, $\alpha = 1.005480$. Yellow curve is residual after the second (final) update of $m$: $m = 0.400499$, $\alpha = 2.010960 $.}
\label{fig:cutted_rn_snr333_res}
\end{figure}

Since $\mu$ and $\eta$ are the same, the truncation lag estimated in Result \ref{thm:result3} (equation \eqref{eqn:lambbd}) is the same, namely $\lambda=0.082$. 
Figures~\ref{fig:cutted_rn_snr333_ew},~\ref{fig:cutted_rn_snr333_ed} and~\ref{fig:cutted_rn_snr333_res} show the truncated wavelets and corresponding predicted data and residuals respectively, at the three iterates of the discrepancy algorithm. Computing the relative data error from the final residual shown in yellow in Figure \ref{fig:cutted_rn_snr333_res}, we conclude that the final slowness estimate $m = 0.400499$ and the final truncated wavelet (yellow curve in Figure \ref{fig:cutted_rn_snr333_ew}) solve the Inverse Problem with $\lambda=0.082, \epsilon=0.27$.

%%%%%%%%%%%
%%%%%%%%%%%

\section{Discussion}
The inverse problem discussed here is a drastic oversimplification of
those encountered in various branches of seismology, as we have
already pointed out. However it is a {\em subproblem} of many active
source inverse problems: limit the data to a single source and a
single receiver, mandate that the material model be spatially
homogeneous, and a very similar problem emerges. Consequently, two
major implication for more complex and prototypical inverse problems
emerge from the theory and examples presented here.

First, variable projection is not enough (to remedy the cycle skipping
misbehaviour of standard FWI). Early contributions to the literature on
VPM applied to seismic inversion concerned linearized or Born
inversion, in which the linear degrees of freedom (analogous to the
source wavelet here) are first order perturbations in material
parameters \cite[]{GockSy:95,vanLeeuwenMulder:09}. This approach appeared in other guises as
Migration-Based Travel Time (MBTT) 
\cite[]{ClementChavent:93} and Reflection FWI (RFWI)
\cite[]{XuChen:12} algorithms.
While these works suggest that
VPM ameliorates the severe nonconvexity in the
nonlinear degrees of freedom (velocity model) noted for basic FWI for reflection data, this mitigation is limited in extent, and cycle-skipping can still occur \cite[]{YinHuang:16}.
VPM reduction of FWI via elimination
of source parameters has been suggested
for example by 
\cite{Rickett:SEG12,Aravkin:12,LiRickettAbubakar:13,Aravkin:GJI16} and
\cite{vanleeuwen2020nonsmooth}, in various contexts.
Result 1 and Figures~\ref{fig:asymmetric_fwi_vpm_01_40hz} and \ref{fig:asymmetric_fwi_vpm_05_40hz} show that for inversion
based on transmitted data, VPM reduction via elimination of source
parameters may fail to produce an optimization problem more tractable
than basic FWI.

Second, the explicit expression \ref{eqn:dexpjgen} for the gradient of $\tJa$ is actually a special case of a generic calculation valid for similar reduced objective functions, based on a large class
of forward maps $F$, as explained by \cite{Symes:SEG20} and
\cite{Symes:22a}, Appendix. This class includes models of transmission
through transparent (non-scattering) material models, of which our
$F[m]$ is the simplest possible case, and some models
of single scattering \cite[]{tenKroode:IPTA14,Symes:IPTA14}. These
calculations, and the features of the reduced objective that they
reveal, are responsible for the feasibility of the VPM approach. We refer
the reader to the cited references for the technical details, and to
\cite{HuangNammourSymesDollizal:SEG19} and references cited there for
numerical illustration in more complex settings.

In the discussion of Experiment 5, we mentioned that the constraint on
noise level posed in Result \ref{thm:result2} is not always necessary for good
convergence to an accurate slowness estimate. Result \ref{thm:result2} presumes a
worst-case constructive buildup of noise impact on $\tJa$: while
sharp, the result may be far too pessimistic for generic noisy data. In
particular, Experiment 5 suggests that uniform filtered random noise in large
quantities (greater than 100\%) still permits fairly precise
determination of slowness. A close reading of the analysis underlying
Result  \ref{thm:result2} \cite[]{Symes:22a} gives some hint about why this might be
so, but much remains to be understood about noise statistics and their
interaction with the solution of the inverse problem.

Result 4 motivates the use of the discrepancy algorithm but does not
really explain its performance. \cite{Symes:21a}, Appendix A,
clarifies the properties of the $\alpha$ update part of the algorithm,
but the interaction with the $m$ update remains obscure, and no
convergence result is proven for the algorithm as a whole in that
work. 

The selection of the error range $[e_-, e_+]$ for the discrepancy algorithm has so
far been entirely ad-hoc, both in earlier work
\cite[]{FuSymes2017discrepancy} and in the experiments reported
here. In almost all practical applications of inverse theory, the
actual data signal-to-noise ratio is a priori unknown. What is the
behaviour of the discrepancy algorithm when the upper error bound
$e_+$ is chosen too low, relative to the actual noise level in the
data? Or the lower error bound $e_-$ too high? Is it possible to
recover from erroneous choice of these parameters, or to learn
anything about the actual noise level from the behaviour of the
discrepancy algorithm? These questions remain
to be investigated.

\section{Conclusion}
We introduced a single-trace acoustic transmission inverse problem, that (despite its simplicity) exhibits one of the fundamental pathologies of real-world full waveform inversion (FWI), namely, cycle-skipping: the standard least-squares objective function of FWI tends to have many stationary points, most of them far from any physically relevant model, so that local optimization produces uninformative model estimates.

We have reviewed the theory of extended source inversion (ESI) applied to this problem. This theory suggests that the ESI approach avoids cycle-skipping, and enables the computation of accurate solutions using standard local optimization methods. Our numerical examples illuminate the extent to which the theory
accurately predicts the performance of ESI in concrete instances. Allowing for the inevitable worst-case flavor of such theory, the conformance of theory and example is reasonably good. In particular, stationary points of the ESI objective function lie near the global minimizer of the usual (FWI) least-squares objective function, with error bounded by a multiple of the wavelet width and the RMS noise level in the data. Extended source inversion involves a penalty formulation so is incomplete without some method for selection of the penalty weight parameter. We presented an algorithm based on the discrepancy principle for adjustment of this parameter, and used it in examples to produce a solution of the inverse problem.

\section{Acknowledgements}
This research is partially supported by the sponsors of the UT Dallas ``3D+4D Seismic FWI'' research consortium.

\section{Data availability statement}
The data that support the findings of this study are available upon reasonable request from the authors.
\clearpage

\bibliographystyle{gji}
\bibliography{masterref,local}

\end{document}